\documentclass[10pt]{amsart}
 \usepackage{amsmath}
 \usepackage{amsxtra}       
 \usepackage{amsthm}        
 \usepackage{amssymb}       
 \usepackage{empheq}        
\usepackage{amscd}
 \usepackage{epsfig}
\usepackage{graphicx}
\usepackage{color}

 \usepackage{mathrsfs} 

 \usepackage{pstricks}     
  \usepackage[all]{xy}       
    \SelectTips{cm}{10}    
    \everyxy={<2.5em,0em>:} 
 \usepackage{hyperref}

 \usepackage{fancyhdr}                   
 \usepackage[paperheight=11in,
             paperwidth=8.5in,
             outer=1.2in,
             inner=1.2in,
             bottom=.7in,
             top=.7in,
             includeheadfoot]{geometry}  
 \renewcommand{\theequation}{\thesection.\arabic{equation}}         
 \makeatletter                                                      
    \@addtoreset{equation}{section} 
    \@addtoreset{footnote}{section} 
    {\hspace*{\fill}$\lrcorner$\endgraf\endgroup\end{trivlist}}     
 \def\newprooflikeenvironment#1#2#3#4{
      \newenvironment{#1}[1][]{
          \refstepcounter{equation}                                 
          \begin{proof}[{\rm\csname#4\endcsname{#2~\theequation}
          \@ifnotempty{##1}{\the\thm@notefont \ (##1)}\csname#4\endcsname{.}}]
          \def\qedsymbol{#3}}
         {\end{proof}}}                                             
 \makeatother                                                       
 \newprooflikeenvironment{definition}{Definition}{$\diamond$}{textbf}%
 \newprooflikeenvironment{example}{Example}{$\diamond$}{textbf}     
 \newprooflikeenvironment{remark}{Remark}{$\diamond$}{textbf}       
 \theoremstyle{plain}                                               
 \newtheorem{thm}[equation]{Theorem}                            

 \newtheorem*{lemma*}{Lemma}                                        
 \newtheorem*{theorem*}{Theorem}                                    
 \newtheorem{lemma}[equation]{Lemma}                                
 \newtheorem{corollary}[equation]{Corollary}                        
 \newtheorem{proposition}[equation]{Proposition}                    
 \DeclareFontFamily{OT1}{pzc}{}                                 
 \DeclareFontShape{OT1}{pzc}{m}{it}{<-> s * [1.100] pzcmi7t}{}  
 \DeclareMathAlphabet{\mathpzc}{OT1}{pzc}{m}{it}                

 \newcommand{\bbar}[1]{\setbox0=\hbox{$#1$}\dimen0=.2\ht0 \kern\dimen0 \overline{\kern-\dimen0 #1}}

 \newcommand{\D}{\mathcal{D}}
 
 \DeclareMathOperator{\End}{\ensuremath{\mathcal{E}\kern-.125em\mathpzc{nd}}}

 \DeclareMathOperator{\Hom}{\mathcal{H}\kern-.125em\mathpzc{om}}

 \newcommand{\I}{\mathcal{I}}

 \newcommand{\m}{\mathfrak{m}}

 \renewcommand{\O}{\mathcal{O}}

 \DeclareMathOperator{\proj}{Proj}
 \DeclareMathOperator{\Proj}{\mathcal{P}\kern-.125em\mathpzc{roj}}

 \renewcommand{\setminus}{\smallsetminus}
 
 \DeclareMathOperator{\spec}{Spec}

 \newcommand{\T}{\mathcal{T}}

 \newcommand{\X}{\mathcal{X}}
 \newcommand{\Y}{\mathcal{Y}}
 \newcommand{\Z}{\mathcal{Z}}

 \newcommand{\cc}{\mathbb{C}}
 \newcommand{\pp}{\mathbb{P}}
  
 \newcommand{\ra}{\rightarrow}

\newcommand{\discr}[3]{\ensuremath{\mathpzc{Discrep}_{#1}^{#3}(#2)}}

 \renewcommand{\th}{\text{th}}
 
 \def\0{\mathbf{0}}
 \def\al{\alpha}
 \def\bet{\beta}							
 \def\del{\delta}
 \def\lam{\lambda}
 \def\lda{\lambda}
 \DeclareMathOperator{\Dim}{dim}
 
 \def\co{\colon\thinspace} 
 \DeclareMathOperator{\Def}{Def}
 \DeclareMathOperator{\card}{card}
  \DeclareMathOperator{\irr}{irr}

 \newcommand{\Mg}[1]{\overline{M}_{#1}}

\author{Maksym Fedorchuk} 
\address{Department of Mathematics \\ Harvard University \\ One Oxford Street \\ Cambridge, MA 02138}
\email{fedorchuk@math.harvard.edu}
\urladdr{\url{http://abel.math.harvard.edu/~fedorchuk}}

\begin{document}
\title{Linear sections of the Severi variety and moduli of curves}

\begin{abstract}
We study the Severi variety $V_{d,g}$ of plane curves of degree $d$ and geometric genus $g$. Corresponding to every such variety, there is a one-parameter family of genus $g$ stable curves whose numerical invariants we compute. Building on the work of Caporaso and Harris, we derive a recursive formula for the degrees of the Hodge bundle on the families in question. For $d$ large enough, these families induce moving curves in $\Mg{g}$.  We use this to derive lower bounds for the slopes of effective divisors on $\Mg{g}$. 
Another application of our results is to various enumerative problems on $V_{d,g}$. 
\end{abstract}
\maketitle

\section{Introduction}

\subsection{Statement of the problem} Let $\pp(d)\cong \pp^{\binom{d+2}{2}-1}$ be the space of all plane curves of degree $d$. Inside it, there is a locally closed subset parameterizing nodal curves with $\delta$ nodes. Its closure $V^{d,\delta}$ is called {\it Severi variety} and has been studied extensively. 
A generic curve in $V^{d,\delta}$ has geometric genus $g=\binom{d-1}{2}-\delta$ and, occasionally, we use $V_{d,g}$ instead of $V^{d,\delta}$. Note that we do not require curves in $V^{d,\delta}$ to be irreducible; the closure of the locus of irreducible curves is denoted by $V^{d,\delta}_{\irr}$. 
Following \cite{CH1}, we denote degrees of $V^{d,\delta}$ and $V^{d,\delta}_{\irr}$ inside $\pp(d)$ by
$N^{d,\delta}$ and $N^{d,\delta}_{\irr}$, respectively.

In this paper, we will be concerned with one-dimensional linear sections of $V^{d,\delta}$ by hyperplanes of a special form in $\pp(d)$. Namely, consider the locus of curves passing through a fixed point $p\in \pp^2$. It is a hyperplane in $\pp(d)$, and is denoted by $H_p$. 

\begin{definition}\label{C-d-delta}
Let $N=\dim V_{d,g}=3d+g-1$ and $\{p_i\}_{1\leq i\leq N-1}$ be a set of general points in $\pp^2$, then
\begin{align*}
C^{d,\delta} &:=V^{d,\delta}\cap H_{p_1} \cap \dots \cap H_{p_{N-1}}, \\
C^{d,\delta}_{\irr} &:=V^{d,\delta}_{\irr}\cap H_{p_1} \cap \dots \cap H_{p_{N-1}}.
\end{align*}
\end{definition}

Let $\Y^{d,\delta}_{\irr}$ be the restriction of the universal family to $C^{d,\delta}_{\irr}$. Then its normalization $(\Y^{d,\delta}_{\irr})^\nu$ is a family of connected, generically smooth curves of genus $g=\binom{d-1}{2}-\delta$, and so induces a rational map from $C_{\irr}^{d,\delta}$ into $\Mg{g}$. The natural question to ask in this situation is, what are the intersection numbers of this curve with the generators of the Picard group of $\Mg{g}$: $$\lambda,\Delta_{0},\Delta_1, \dots, \Delta_{\lfloor g/2\rfloor}?$$

After a moment of reflection, we see that intersection numbers of $C^{d,\delta}_{\irr}$ with boundary divisors can be expressed in terms of the degrees of various Severi varieties. For example, 
$$C^{d,\delta}_{\irr}\cdot \Delta_0=(\delta+1)N^{d,\delta+1}_{\irr}.$$
The more subtle problem is determining the degree of the Hodge bundle on $C^{d,\delta}_{\irr}$, i.e., the number $$C^{d,\delta}_{\irr} \cdot \lambda,$$ which we denote by $L^{d,\delta}_{\irr}$.

This question is even more intriguing when the Brill-Noether number $\rho(g,2,d)=3d-2g-6$ is non-negative, in which case the image of $C^{d,\delta}_{\irr}$ is a moving curve inside $\Mg{g}$, and so the fraction
\begin{align*}
\frac{C^{d,\delta}_{\irr} \cdot \Delta}{L_{\irr}^{d,\delta}}
\end{align*}
gives a lower bound on the slope of effective divisors on $\Mg{g}$. 
 
 The main result of this paper is the recursive formula of Theorem \ref{main theorem}  that allows us to compute numbers $L^{d,\delta}_{\irr}$. The statement of the theorem requires several preliminary definitions which we give in the remainder of this section. 
 
 \subsection{Notations and conventions} We work over the field of complex numbers $\cc$. We denote by $\Mg{g,n}$ the coarse moduli space of stable curves of genus $g$ with $n$ marked points. For $g\geq 2$, we let $\mathcal{C}_g \stackrel{\pi}{\rightarrow} \Mg{g}$ be the ``universal'' curve, defined away from the codimension two locus of curves with extra automorphisms. Then the Hodge bundle is, by definition, $\mathbb{E}:=\pi_*(\omega_{\mathcal{C}_g/\Mg{g}})$. We set $\lambda:=c_1(\mathbb{E})$.
 
 Throughout the paper, $\delta_{ij}$ stands for the Kronecker's delta. We reserve the symbol $\epsilon$ for the generator of the ring of dual numbers $\spec \cc[\epsilon]/(\epsilon^2)$. Given a scheme $X$, we denote its normalization by $X^{\nu}$. The tangent space to $X$ at a point $x$ is denoted by $\mathbb{T}_x X$.
 The tacnode of order $m$ is a planar curve singularity analytically isomorphic to the singularity of $y^2-x^{2m}=0$ at the origin. For the sake of uniformity, we will not distinguish between a node and a tacnode of order $1$. 
 \subsection{Reducible vs. irreducible curves}\label{red-vs-irred}
 
 For technical reasons, it is simpler to work with the variety $V^{d,\delta}$ of possibly reducible curves. 
 Recall from Definition \ref{C-d-delta}, that the curve $C^{d,\delta}$ parameterizes $\delta$\nobreakdash-nodal curves of degree $d$ through $N-1$ general points $\{p_i\}_{1\leq i \leq N-1}$ in $\pp^2$. Consider a component of $C^{d,\delta}$ parameterizing curves which have $k$ irreducible components belonging to Severi varieties $V_{\irr}^{d_1,\delta_1}, \dots, V_{\irr}^{d_k, \delta_k}$ where  
 \begin{equation}\label{degrees}
 \begin{split}
 d &=\sum_{i=1}^k d_i \ \ \text{and} \\
 \delta &=\sum_{i=1}^k \delta_i +\sum_{1\leq i <j \leq k} d_id_j. 
 \end{split}
 \end{equation}

 For $1\leq i\leq k$, we set $g_i=\binom{d_i-1}{2}-\delta_i$ and $N_i=\dim(V_{d_i,g_i})=3d_i+g_i-1$. Equalities \eqref{degrees} imply that $\sum_{i=1}^k g_i -k+1=g$. Clearly, the component in question is a union of the Segre images of the products 
 $$\prod_{i=1}^{k} (V_{\irr}^{d_i,\delta_i}\cap (\cap_{j=1}^{N_i-\delta_{i,r}} H_{p_{i,j}})) = \{pts\} \times \cdots \times C_{\irr}^{d_r,\delta_r} \times \cdots \times \{pts\} \subset \prod_{i=1}^{k} V_{\irr}^{d_i,\delta_i} ,$$
 for every $1\leq r\leq k$, and where $\{p_{i,j}\} 
 =\{p_1,\dots,p_{N-1}\}$ (we use the shorthand $\{pts\}$ to denote a finite set of isolated points on a variety).
 
 There is a natural map 
 $$j: \prod_{i=1}^{k} V_{\irr}^{d_i,\delta_i} \dasharrow \prod_{i=1}^k \Mg{g_i}$$
 defined as the product of moduli maps.  Define a line bundle $\Lambda$ on $\prod_{i=1}^k \Mg{g_i}$ to be the product of pullbacks of the determinants of the Hodge bundles from the individual factors.
 
  Define now an auxiliary space $\mathcal{M}^g$ to be a countable union of varieties $\prod_{i=1}^k \Mg{g_i}$ over all $k$ and $\{g_i\}_{1\leq i\leq k}$ satisfying $\sum_{i=1}^k g_i-k+1=g$. As described in the previous paragraph,  $\mathcal{M}^g$ comes with the line bundle $\Lambda$ whose first Chern class we denote by $\lambda$. By above we have a natural ``moduli" map $$j:C^{d,\delta} \dasharrow \mathcal{M}^g.$$
  Note that for a given pair $(d,\delta)$ we need only a finite number of components of $\mathcal{M}^g$.
  Therefore, we can define the intersection number
  $$L^{d,\delta}:=j_*(C^{d,\delta})\cdot \lambda.$$
 Using this definition, we have
 $$j_*(\{pt\} \times \cdots \times C^{d_r,\delta_r} \times \cdots \times \{pt\})\cdot \lambda=L_{\irr}^{d_r,\delta_r},$$ where $L_{\irr}^{d_r,\delta_r}$ is defined via the map 
 $$j_r: C_{\irr}^{d_r,\delta_r} \dasharrow \Mg{g_r}.$$ 
 
Numbers $L^{d,\delta}_{\irr}$ are recovered inductively from numbers $L^{d,\delta}$ using the following formula
 \begin{align*}
 L^{d,\delta}_{\irr} &= L^{d,\delta} - \sum_{\substack{ k \\ (d_1,\dots,d_k) \\ (\delta_1, \dots, \delta_k) } } \sum_{r=1}^k \left(\prod_{\substack{i=1, \\ i\neq r}}^n \binom{3d+g-1}{3d_i+g_i-1} N_{\irr}^{d_i,\delta_i} \right) \binom{3d+g-1}{3d_r+g_r-2} L_{\irr}^{d_r,\delta_r} \ ,
 \end{align*}
  where the sum is taken over all $k$ and over all ordered partitions $(d_1,\dots, d_k)$  of $d$, and $(\delta_1, \dots, \delta_k)$ of $\delta-\sum_{1\leq i < j \leq k} d_id_j$, respectively; and 
  where $g_i=\binom{d_i-1}{2}-\delta_i$.
 
\subsection{Generalized Severi variety} We base our analysis of $C^{d,\delta}$ on the degeneration approach developed by Caporaso and Harris in \cite{CH1}. First, we recall the definition of the generalized Severi variety from ibid. Section 1.1, and the notations accompanying it. 

Given a sequence of non-negative integers  $\alpha=(\alpha_1,\alpha_2, \dots)$, we define 
\begin{align*}   | \alpha | =\sum_i \alpha_i, &&  I\alpha=\sum_i i\alpha_i, && I^{\alpha}&=\prod_i i^{\alpha_i}.
\end{align*}
Fx a line $L\subset \pp^2$, once and for all. 

\begin{definition} For a given $d$ and $\delta$, 
consider any two sequences $\alpha=(\alpha_1,\alpha_2, \dots)$ and $\beta=(\bet_1,\bet_2,\dots)$ of non-negative integers such that 
\begin{align*}
I\alpha+I\beta=d.
\end{align*} 
Fix a general collection of points 
\begin{align*}
\Omega=\{ p_{i,j}\}_{1\leq j\leq \alpha_j} \subset L .
\end{align*}
Consider the locus of $\delta$-nodal plane curves $X$ of degree $d$ that does not contain $L$, and such that, for the normalization map 
$$\eta: X^{\nu}\ra X,$$ we have 
$$\eta^*(L)=\sum i\cdot q_{i,j}+\sum i\cdot r_{i,j}$$ 
for some $|\al|$ points $q_{i,j}$ and $|\bet|$ points $r_{i,j}$ on $X^{\nu}$ such that 
$$\eta(q_{i,j})=p_{i,j}.$$
The closure of this locus is called the {\it generalized Severi variety} and is denoted 
\begin{align*}
V^{d,\delta}(\al,\bet)(\Omega).
\end{align*}
\end{definition}

\begin{remark}
Since $\Omega$ is a general set, the geometry of $V^{d,\delta}(\al,\bet)(\Omega)$ does not depend on it. Therefore, it is customary to omit $\Omega$ from the notation.
\end{remark}
\begin{definition}
We also define 
\begin{align*}
V^{d,\delta}_L(\al,\beta):=\{X\cup L \ : \ X\in V^{d,\delta}(\al,\beta) \} \subset \pp(d+1). 
\end{align*}
\end{definition}

We recall next the main result of \cite{CH1}:
\begin{thm}\cite[Theorem 1.2]{CH1}
\label{CH hyperplane section theorem}
For a general $q\in L$, we have the following equality of cycles
\begin{equation}\label{cycles}
\begin{split}
V^{d,\delta}(\alpha,\beta)(\Omega) \cap H_q&= \sum_k kV^{d,\delta}(\alpha+e_k,\beta-e_k)(\Omega \cup \{q\}) \\
&+\sum I^{\beta'-\beta}\binom{\beta'}{\beta}V_L^{d-1,\delta'}(\alpha',\beta')(\Omega'),
\end{split}
\end{equation}
where the second sum is taken over all triples $(\delta',\alpha',\beta')$ satisfying \footnote{It follows from Equation \eqref{cycles} that only varieties $V_L^{d-1,\delta'}(\alpha',\beta')(\Omega')$ with $\alpha'\leq \alpha, \beta'\geq \beta$ and $\delta'\leq \delta$ appear with the non-zero coefficient. }
\begin{align*} 
  &|\beta'-\beta|+\delta-\delta'=d-1,
\end{align*} and over all sets of points $\Omega'=\{ p'_{i,j}\}_{1\leq j \leq \alpha'_i} \subset \Omega$. 
\end{thm}
We refer to components of $H_q:=V^{d,\delta}(\alpha,\beta)(\Omega) \cap H_q$ appearing on the first line of Equation \eqref{cycles} as {\it Type I components}. A generic point of a Type I component does not contain the line $L$. The remaining components of $H_q$ are called {\it Type II components} and parameterize curves containing the line $L$. 

An immediate corollary of Theorem \ref{CH hyperplane section theorem} is 
\begin{thm}\cite[Theorem 1.1]{CH1}
\label{CH recursive formula}
The degrees $N^{d,\delta}(\alpha,\bet)$ of the Severi varieties $V^{d,\delta}(\alpha,\beta)$ satisfy the recursion
\begin{align*} 
N^{d,\delta}(\alpha,\beta) &= \sum_k kN^{d,\delta}(\alpha+e_k,\beta-e_k) \\
&+\sum I^{\beta'-\beta}\binom{\alpha}{\alpha'}\binom{\beta'}{\beta}N^{d-1,\delta'}(\alpha',\beta'),
\end{align*}
where the second sum is taken over all triples $(\delta',\alpha',\beta')$ satisfying $|\beta'-\beta|+\delta-\delta'=d-1.$
\end{thm}

Throughout the paper we work with linear sections of $V^{d,\delta}(\alpha,\beta)$. We introduce the following notations.
\begin{definition}\label{C-d-delta-2}
Let $N:=\dim V^{d,\delta}(\alpha,\beta)$ and $\{p_i\}_{1\leq i\leq N-1}$ be the set of $N-1$ generic points of $\pp^2$. Set
\begin{align*}
S^{d,\delta}(\alpha,\beta) &:=V^{d,\delta}(\alpha,\beta)\cap H_{p_1} \cap \dots \cap H_{p_{N-2}} \notag\\
\intertext{and}
C^{d,\delta}(\alpha,\beta) &:=V^{d,\delta}(\alpha,\beta)\cap H_{p_1} \cap \dots \cap H_{p_{N-1}}. \notag
\end{align*}
\end{definition}
Note that $S^{d,\delta}(\alpha,\beta)$ is a surface and $C^{d,\delta}(\alpha,\beta)$ is a divisor on it.
\begin{definition}
For $g=\binom{d-1}{2}-\delta$, we let $j$ to be the induced rational map 
$$j:C^{d,\delta}(\alpha,\beta) \dasharrow \mathcal{M}^g,$$
where $\mathcal{M}^g$ is the scheme described in Section \ref{red-vs-irred}. 
We  set
$$L^{d,\delta}(\alpha,\beta):=j_*(C^{d,\delta}(\alpha,\beta))\cdot \lambda.$$
\end{definition}

The strategy for calculating numbers $L^{d,\delta}$ is now clear. Using the degeneration of Theorem \ref{CH hyperplane section theorem}, we produce a recurrence relation among numbers $L^{d,\delta}(\alpha,\beta)$ paralleling that of Theorem \ref{CH recursive formula}. The main theorem of this paper is
\begin{thm}\label{main theorem} The numbers $L^{d,\delta}(\alpha,\beta)$ satisfy the recursion
\begin{align*} 
L^{d,\delta}(\alpha,\beta) &= \sum_k kL^{d,\delta}(\alpha+e_k,\beta-e_k) \\
&+\sum I^{\beta'-\beta}\binom{\alpha}{\alpha'}\binom{\beta'}{\beta}L^{d-1,\delta'}(\alpha',\beta')  \\
&+\frac{1}{12}\sum I^{\beta'-\beta}\binom{\alpha}{\alpha'}\binom{\beta'}{\beta}\cdot\biggl(\sum_{k} (\beta'_{k}-\beta_{k})(k^2-1)\biggr)\cdot N^{d-1,\delta'}(\alpha',\beta') \ ;\end{align*}

where the second sum is taken over all triples $(\delta',\alpha',\beta')$ satisfying $|\beta'-\beta|+\delta-\delta'=d-1$ and the third sum is taken over all triples $(\delta', \alpha', \beta')$ satisfying
$|\beta'-\beta|+\delta-\delta'=d-2$.
\end{thm}

Intuitively, the first and the second line of the recursion do not require an explanation. Given Theorem \ref{CH hyperplane section theorem}, they are, at least, expected to appear. The third line of the recursion is of different nature. It arises from an extra complication that did not appear in the analysis of \cite{CH1}. We pause here to describe it.

The degree is an intrinsic property of the Severi variety, as it comes naturally with the definition. On the other hand, the first Chern class of the Hodge bundle on $C^{d,\delta}(\al,\bet)$ is defined in terms of an extra structure that we put on the Severi variety, 
namely, the moduli map to $\Mg{g}$. Moreover, the moduli map is not rational. It would not pose much difficulty if we were considering $C^{d,\delta}(\al,\bet)$ alone, as the map would then naturally extend to the regular map (at least after the normalization). However, our approach is to degenerate $C^{d,\delta}(\al,\bet)$ to a union of linear sections of ``simpler'' generalized Severi varieties. The total space $S$ of the degeneration is thus two-dimensional, and we no longer can expect the moduli map to extend. To resolve the moduli map, we have to blow-up $S$. The exceptional divisors of the blow-up will then contribute to the calculation. These contributions appear on the third line of the recursion in Theorem \ref{main theorem}.

Yet another way to think of the problem is to recall that we have the Kontsevich moduli space of stable maps $\Mg{g,0}(\pp^2,d)$ that fits into the diagram
\begin{align}
 \xymatrix {\Mg{g,0}(\pp^2,d) \ar[d]^{\pi}  \ar[r] & \Mg{g} \\ 
 V_{d,g} \ar@{-->}[ur]  
 }
 \end{align}
The projection $\pi$ maps the so called ``main component'' of $\Mg{g,0}(\pp^2,d)$ birationally onto $V_{d,g}$.
The space of stable maps has the advantage that the moduli map to $\Mg{g}$ is a well-defined morphism. The difficulty that arises, if one wants to work with the stable maps, is the existence of components of $\Mg{g,0}(\pp^2,d)$ parameterizing maps that contract components of positive genus. Moreover, these components of $\Mg{g,0}(\pp^2,d)$ have wrong dimension. 

The moduli map from $S$ to $\Mg{g}$ will be undefined precisely at the points corresponding to the stable maps contracting components of positive genus. To resolve the map, we are thus required to understand the proper transform of $S$ inside $\Mg{g,0}(\pp^2,d)$. 

Finally, even though we do not use the language of stable maps in this paper, the above discussion serves as a motivation for much of what follows. 

\subsection{Structure of the paper:} In Chapter \ref{remarks}, we explain the notions of $\lambda$-indeterminacy and discrepancy that we encounter when working with two-dimensional families of not necessarily stable curves. In Chapter \ref{global geometry}, we describe explicitly where the indeterminacy occurs along a special hyperplane section of $S^{d,\delta}(\alpha,\beta)$. In Chapter \ref{local geometry}, we reduce the discrepancy calculation on $S^{d,\delta}(\alpha,\beta)$ to a discrepancy calculation on a surface in the product of deformation spaces of several tacnodes. We develop a theory which allows to do this local calculation in Chapters \ref{single tacnode} and \ref{multiple tacnodes}, and perform the calculation in Chapter \ref{calculation}. In the final chapter of the paper, we give examples of enumerative problems that can be solved using our Theorem \ref{main theorem}.

\subsection{Acknowledgments:} We are grateful to Joe Harris for numerous discussions, advice, and for introducing us to this problem. We we would like to thank Ethan Cotterill, for reading the preliminary draft and providing valuable suggestions, and Anatoly Preygel, whose Python script served as a basis for the program\footnote{The code is available upon request.} implementing the recursion of
our Theorem \ref{main theorem}. 

\section{Remarks on the intersection theory}\label{remarks}
In this chapter, we recall some generalities on the coarse moduli space $\Mg{g}$ of stable curves of genus $g$. The reference for the material presented here is \cite[Chapter 3]{HM-book}. 

\subsection{Families of curves over one-dimensional bases.} 
Consider a proper family $\X \ra C$ of curves over a one-dimensional irreducible base $C$. Suppose that a generic fiber is a nodal curve of arithmetic genus $g$. Then there is a finite surjective base change $f: C' \ra C$, and a family $\pi: \Y\ra C'$ of stable curves of genus $g$ such that $\Y$ is birational to $\X\times_C C'$. Colloquially, $\Y\ra C'$ is a stable reduction of $\X\ra C$. The Hodge bundle $E_\Y$ on $C'$ is defined by
$$E_\Y:=\pi_*(\omega_{\Y/C'}).$$
It is a pullback of the Hodge bundle on $\Mg{g}$ under the natural morphism 
$C'\ra \Mg{g}$ induced by the family $\Y\ra C'$. Without performing a stable reduction, we still have a natural rational map\footnote{When $C$ is smooth, the map extends to a regular morphism.} from $C$ to  $\Mg{g}$:
$$j: C\dasharrow \Mg{g}.$$

\begin{definition} The {\bf degree of $\lambda$ on $\X/C$}, denoted by $\lambda_{\X/C}$, is the intersection number $j_*(C)\cdot \lambda$. It is also equal to $c_1(E_\Y)/\deg f$.
\end{definition}
Note that $\lambda_{\X/C}$ depends only on the geometry of the family at a generic point of $C$. If the family $\X\ra C$ is understood, we use the shorthand $\lambda_{C}$ to denote the degree of $\lambda$ on $\X/C$.

The following two lemmas follow from definitions.
\begin{lemma}\label{lambda-union}
Suppose $\X\ra C$ is a flat proper family of curves with a nodal generic fiber. If $\X=\X_1\cup \X_2$, a union of two families, then 
$$\lambda_{\X/C}= \lambda_{\X_1/C}+\lambda_{\X_2/C}.$$
\end{lemma}

\begin{lemma}\label{lambda-nodes}
Suppose $\X\ra C$ is a flat proper family of curves of arithmetic genus $g$. Suppose, moreover, that a generic fiber of $\X$ has $\delta$ nodes and no other singularities. Then the normalization $\X^{\nu}\ra C$ is a proper family, whose generic fiber is a smooth curve of genus $g-\delta$, and we have
$$\lambda_{\X/C}=\lambda_{\X^\nu/C}.$$ 
\end{lemma}

\subsection{Families of curves over two-dimensional bases.} 

Consider now a flat proper family $\X \ra B$ over a two-dimensional irreducible base $B$ whose generic fiber is a stable curve of arithmetic genus $g$. We have an induced rational map $$j: B\dasharrow \Mg{g},$$
which we call the {\it moduli map}. The locus $U$ of points $b$ such that $\X_b$ is stable is an open subset of $B$. Note that $j$ is regular on $U$. Replace now $B$ by its normalization, and denote by $\X$ its own pullback to the normalization, and by $U$ its own preimage. 

By properness of $\Mg{g}$, the moduli map $j$ is defined away from a finite set of points of $B$. We call them {\it points of indeterminacy} of the moduli map, and denote the set of such points by $\text{\underline{Indet}}(\X/B)$; or $\text{\underline{Indet}}(B)$ if the family $\X$ is understood. Set $W=B\setminus \text{\underline{Indet}}(B)$. Then $U$ is a subset of $W$, and the inclusion can be proper. For example, suppose $b\in B\setminus U$ is such that the isomorphism class of the stable limit of any one-parameter family with the center at $b$ does not depend on the family. Then $j$ extends to a regular map in a neighborhood of $b$, even though $\X_b$ is not stable. 

The resolution of $j$ is a proper birational map $\pi: \hat{B} \ra B$, restricting to the isomorphism on $U$, together with a regular map $\hat{\jmath}: \hat{B} \ra \Mg{g}$ such that the following diagram commutes:
\begin{align}
 \xymatrix {\hat{B} \ar[d]^{\pi} \ar[dr]^{\hat{\jmath}} \\ 
 B \ar@{-->}[r]^{j} & \Mg{g}
 }\label{c}
\end{align}
By the Zariski's Main Theorem, $\pi: \hat{B}\ra B$ has connected fibers.
 We note that it is possible to test whether $b\in B$ will be a point of indeterminacy without passing to the normalization. Suppose that the isomorphism class of the stable limit of any one-parameter family with the center at $b$ belongs to a discrete set. Then, after the normalization, the moduli map is defined at $b$. Indeed, the finiteness assumption implies that $\pi^{-1}(b)$ is a finite set, and so $\pi$ must be an isomorphism at $b$.

For every curve $C\subset B$, we can define the number $\lambda_C$ by restricting $j$ to $C$ and setting $$\lambda_C:=j_*(C)\cdot \lambda.$$ Note that if $\X\times_B C$ is a family with a nodal generic fiber, then 
$\lambda_C=\lambda_{\X\times_B C /C}$. 

For every $Q$-Cartier divisor $C\subset B$, we define another closely related number.
\begin{definition} The number $$(\lambda\cdot C)_B:=\hat{\jmath}^{\, *}(\lambda)\cdot\pi^*(C)$$ is called the {\it $\lambda$-degree of $C$ on $B$}.
This also equals to $\lambda_{C'}$ for any divisor $C'$ that does not pass through $\text{\underline{Indet}}(B)$ and is linearly equivalent to $C$.
\end{definition}


For a $Q$-Cartier divisor $C$, the pullback is defined and we have 
$$\pi^*(C)=\pi_*^{-1}C+E_C,$$ where $E_C$ is some linear combination of exceptional divisors of $\pi$. We define the $\lambda$-discrepancy, or simply discrepancy, along $C$ to be the number
$$\discr{B}{C}{\lambda}:=(\lambda\cdot C)_B-\lambda_C.$$
For $b \in \text{\underline{Indet}}(B)$, we define $$\discr{b}{C}{\lambda}:=\sum E \cdot \lambda$$ where the sum is taken over all exceptional divisors $E$ mapping to $b$.

Note that since $\Mg{g}$ is only a coarse moduli space, there might not be a family of stable curves over $\hat{B}$ that induces the map $\hat{\jmath}$. However, by \cite[Lemma 3.89]{HM-book}, there is a finite order base change $\hat{\hat{B}} \xrightarrow{f} \hat{B}$ and a family of stable curves $\Y \ra \hat{\hat{B}}$ such that 
$\hat{\jmath}\circ f: \hat{\hat{B}} \ra \Mg{g}$ is induced by the family $\Y$. Moreover, we have $$\X\times_U f^{-1}(U) \cong \Y_{f^{-1}(U)}.$$ 
We set $\hat{\pi}:=\pi\circ f$.
By abuse of terminology, we call any such family, $\Y\ra \hat{\hat{B}}$, a {\it stable reduction} of $\X\ra B$ (cf. \cite[Corollary 3.96]{HM-book}).

Given a Q-Cartier divisor $C\subset B$, we define its strict transform on $\hat{\hat{B}}$ to be $\overline{C}:=\overline{\hat{\pi}^{-1}(C\cap W)}$. We define the exceptional part of $C$ to be $\hat{E}_C:=\hat{\pi}^{\, *}(C)-\overline{C}.$ Then the $\lambda$-discrepancy along $C$ on $B$ equals to $$\frac{1}{\deg{f}}(\hat{E}_C\cdot \lambda).$$
Note that $\hat{E}_C\cdot \lambda$ is defined unambiguously. It is the degree of the Hodge bundle corresponding to the family of stable curves, $\Y_{\hat{E}_C}\ra \hat{E}_C$.

\begin{lemma}\label{indeterminacy}
If $b\in B$ is such that the geometric genus of the fiber $\X_b$ is $g$, then $b \notin \text{\underline{Indet}}(B)$.
\end{lemma}

\begin{proof}
Suppose $\pi: \hat{B}\ra B$ is a minimal resolution of the moduli map, and $E\subset \hat{B}$ is an exceptional divisor mapping to $b$. After a finite base change, we can assume that we have a family $\Y$ of stable curves over   
$\hat{B}$. The fibers of $\Y$ over $E$ are stable curves of arithmetic genus $g$ that map to $\X_b$. Since the normalization $(\X_b)^{\nu}$ of $\X_b$ has geometric genus $g$, we conclude that all fibers of $\Y$ over $E$ must be isomorphic to $(\X_b)^{\nu}$. Hence, $E$ maps to a point in $\Mg{g}$ and so cannot be an exceptional divisor of $\pi$. 
\end{proof}

Given two linearly equivalent divisors $C$ and $C'$ on $B$, we have, by definition, $(\lambda\cdot C)_B=(\lambda\cdot C')_B$. If, moreover, $C'$ lies in $U$, we can compute $\lambda_C$ in terms of $\lambda_{C'}$ and the discrepancy $\discr{S}{C}{\lda}$. This simple observation is important because the discrepancy along $C$ depends only on the geometry of the family $\X\ra B$ in the analytic neighborhood of $B$ around $\text{\underline{Indet}}(B)$. Therefore, to compute the discrepancy we can work locally around each point of indeterminacy. 




\section{The degeneration and points of indeterminacy}\label{global geometry}

Consider the surface
\begin{align*}
S^{d,\delta}(\alpha,\beta) =V^{d,\delta}(\alpha,\beta)\cap H_{p_1} \cap \dots \cap H_{p_{N-2}} 
\end{align*}
and a restriction of the universal curve to $S^{d,\delta}(\al,\beta)$ which we denote by $\Y:=\Y^{d,\delta}(\al,\beta)$. A generic fiber of $\Y$ is a $\delta$-nodal curve of geometric genus $g=\binom{d-1}{2}-\delta$. By discussion in Chapter \ref{remarks}, we have an induced rational moduli map $S^{d,\delta}(\alpha,\beta) \dasharrow \mathcal{M}^{g}$.

Let $p_{N-1}$ be a general point in $\pp^2$. Then by Definition \ref{C-d-delta-2}, 
\begin{align}
\label{gensec}
S^{d,\delta}(\alpha,\beta)\cap H_{p_{N-1}}= C^{d,\delta}(\alpha,\beta).
\end{align}

Suppose now $q$ is a general point on the line $L$. Then by the Theorem \ref{CH hyperplane section theorem}, we have
\begin{equation}\label{spsec}
\begin{split}
S^{d,\delta}(\alpha,\beta)\cap H_q &=\sum_k kC^{d,\delta}(\alpha+e_k,\beta-e_k) \\
&+\sum I^{\beta'-\beta}\binom{\alpha}{\alpha'}\binom{\beta'}{\beta}C_{L}^{d-1,\delta'}(\alpha',\beta'). 
\end{split}
\end{equation}

Equalities \eqref{gensec} and \eqref{spsec} show that there is a linear equivalence 
\begin{equation}\label{linear-equiv}
\begin{split}
C^{d,\delta}(\alpha,\beta)& \sim \sum_k kC^{d,\delta}(\alpha+e_k,\beta-e_k) \\
&+\sum I^{\beta'-\beta}\binom{\alpha}{\alpha'}\binom{\beta'}{\beta}C_{L}^{d-1,\delta'}(\alpha',\beta')
\end{split}
\end{equation}
of divisors on $S^{d,\delta}(\alpha,\beta)$.

By results of \cite{CH1}, the stable limit of a generic arc with the center at a generic point $[C]=[Y\cup L]$ of $C_{L}^{d-1,\delta'}(\alpha',\beta')$ is a nodal curve, 
which is a partial normalization of $C$, equal to a union of the line $L$ and the normalization of the curve $X$. 
Recalling the discussion of Chapter \ref{remarks}, we have the following result.
\begin{corollary}\label{indet-type2} After the normalization of $S^{d,\delta}(\al,\bet)$, the moduli map is defined at a generic point of any Type II component.
\end{corollary}

Invoking Lemmas \ref{lambda-union} and \ref{lambda-nodes}, we have
\begin{align*}
\lambda_{C_{L}^{d-1,\delta'}(\alpha',\beta')}=L^{d-1,\delta'}(\alpha',\beta').
\end{align*}

Finally, we have, by the definition, 
\begin{equation}\label{L-induction}
\begin{split} 
&L^{d,\delta}(\alpha,\beta)+\discr{S^{d,\delta}(\alpha,\beta)}{C^{d,\delta}(\al,\bet)}{\lambda} \\ &= \sum_k kL^{d,\delta}(\alpha+e_k,\beta-e_k) 
+\sum I^{\beta'-\beta}\binom{\alpha}{\alpha'}\binom{\beta'}{\beta}L^{d-1,\delta'}(\alpha',\beta') 
+\discr{S^{d,\delta}(\alpha,\beta)}{H_q}{\lambda}.
\end{split}
\end{equation}
Therefore to prove Theorem \ref{main theorem}, we need to understand the points of indeterminacy $\text{\underline{Indet}}(S^{d,\delta}(\alpha,\beta))$ of the moduli map $j$ and the $\lambda$-discrepancies along $C^{d,\delta}(\al,\bet)$ and $H_q$.

\subsection{Indeterminacy Points}\label{indeterminacy points}

In this section we describe the points of indeterminacy that occur along linear sections $H_{p_{N-1}}$ and $H_{q}$.  The main tools in our analysis are nodal reduction for curves and the following fundamental dimension-theoretic result on the deformations of plane curves:

\begin{lemma}
\label{dimLemma}
\cite[Corollary 2.7]{CH1}
Fix a subset $\Omega$ of general points on $L\subset \pp^2$. Let $V$ be a locally closed irreducible subvariety of $\pp(d)$ and $X$ a generic point of $V$. Let $\pi: X^{\nu}\ra X$ be the normalization map and
$e:=\card(X\cap (L\setminus \Omega))$. Then
\begin{align}\label{dim-ineq}
\Dim V \leq 2d+g-1+e.
\end{align}
Moreover, if equality holds and $\card( \pi^{-1}(L\setminus\Omega))=e$, then $V$ is a dense open subset of a generalized Severi variety. 
\end{lemma}

The following is the main result of this chapter. 
\begin{proposition}\label{points of indeterminacy}
The generic hyperplane section $C^{d,\delta}(\al,\beta)$ does not contain points of indeterminacy. The points of indeterminacy that lie on $H_{q}$ are $$V^{d,\del'}(\alpha',\beta')\cap H_{p_1}\cap\dots\cap H_{p_{N-2}},$$
where $(\del', \alpha', \beta')$ satisfy 
\begin{align}
\label{indeterminacy condition}
|\beta'-\beta|+\del-\del_0=d-2, \text{and} \ \beta'_i>0 \ \text{for some}\  i\geq 2.
\end{align}

 \end{proposition}
 
\begin{remark}
The points of indeterminacy on $H_q$ can also be described as points of $$V^{d,\delta'}(\alpha',\beta'-e_i-e_j+e_{i+j})\cap H_{p_1}\cap\dots\cap H_{p_{N-2}},$$
where $V^{d,\delta'}(\al',\bet')$ is a {\it Type II} component of $H_{q}$.
\end{remark}

\subsection{ Nodal reduction:}\label{stable reduction} In our analysis of indeterminacy points on $H_{q}$, we will often need to understand the stable limits of one-parameter families inside $S^{d,\delta}(\al,\bet)$. We use the approach of Section 3 of \cite{CH1} in what follows. Let $\Delta_0$ be an arc in $S^{d,\delta}(\al,\bet)$ with the center at $[C]=[Y\cup L] \in C_L^{d-1,\del'}(\al',\bet')$ and $\Y=\Y^{d,\delta}(\al,\bet)\times_{S} \Delta_0$ the restriction of the universal family.

It is a standard result that we may perform a nodal reduction of $\Y\ra \Delta_0$ to obtain a flat family $\Z \ra \Delta$ of generically smooth genus $g$ curves satisfying the following conditions:
\begin{enumerate}
\item The total space $\Z$ is smooth and there is a map
 $$\eta: \Z \ra \pp^2.$$
\item The central fiber $Z_0$ is nodal.
\item $Z_0$ decomposes as $$Z_0=\overline{Y}\cup P,$$ where $\overline{Y}$ is a strict transform of $Y$, and $P$ is a union of components mapping to $L$.

\item There is a multisection $F=\sum Q_{i,j}$ such that $\eta(Q_{i,j})=p_{i,j}$, and a multisection $V= \sum R_{i,j}$ such that
$$\eta^{*}(L)=\sum_{i,j} i\cdot Q_{i,j} + \sum_{i,j} i\cdot R_{i,j} + P_0,$$
where $P_0$ is a divisor supported on $P$. 

\end{enumerate}
\[
\xymatrix{
Z_0\ar@<-2ex>[d] \subset \Z \ar[r]^{\eta} \ar@<3ex>[d] &\pp^2 \\
 [C] \in \Delta & 
}
\]

We consider a decomposition 
$$P=\sum_{i=1}^{c_1} P'_i +\sum_{i=1}^{c_2} P''_i$$ 
of $P$ into connected components such that $P'_i$'s are not contracted by $\eta$ and $P''_i$'s are contracted by $\eta$.

Every contracted component $P''_i$ has to meet either $V$ or $F$. We let the union of these to be $P''_v$ and $P''_f$, respectively. Denote the number of connected components of $P''_v$ and $P''_f$ by $v$ and $f$, respectively.

\hfill

{\it Proof of Proposition \ref{points of indeterminacy}:}

\hfill

First, we consider $H_{p_{N-1}}$. Note that $H_{p_{N-1}}$ intersects any Type II component at a generic point of the component. By Corollary \ref{indet-type2}, the moduli map is defined at these points. All other points of $H_{p_{N-1}}$ correspond to curves not containing the line $L$, and we next consider only such points. By Lemma \ref{dimLemma}, the only points on $H_{p_{N-1}}$ where the geometric genus drops are $(\delta+1)$-nodal curves. The moduli map is clearly defined at these points. By Lemma \ref{indeterminacy} there are no other points of indeterminacy.

\hfill

{\it Indeterminacies on Type I components:} 

Consider a curve $C$ on a Type I component $C^{d,\del}(\al+e_k, \bet-e_k)$ of $S^{d,\del}(\al,\bet) \cap H_q$. By Lemma \ref{indeterminacy}, a point of indeterminacy can occur only when the geometric genus drops. If $C$ does not contain a line, by Lemma \ref{dimLemma}, genus can drop by at most $1$. In this case we must have $\card(C\cap (L\setminus \Omega))=|\beta-e_k|$, and the inequality in \eqref{dim-ineq} becomes equality. Therefore, $C$ is a $(\delta+1)$-nodal curve, smooth along $L$. 
We conclude that potential points of indeterminacy have to contain the line $L$. By Theorem \ref{CH hyperplane section theorem}, these are points in $V^{d-1,\del'}(\al',\bet')$ such that 
$$|\bet'-(\beta-e_k)|+\del-\del'=d-1  \  \Leftrightarrow \ |\bet'-\beta|+\del-\del'=d-2.$$ 

\hfill

{\it Indeterminacies on Type II components:} 

Suppose $[C] \in C_{L}^{d-1,\del'}(\al',\bet')$ is a point of indeterminacy. By definition, $C=Y\cup L$ where $Y$ is in $C^{d-1,\del'}(\al',\bet')$, for some $(\del',\al',\bet')$ satisfying
\begin{align}\label{CND}
|\bet'-\bet|+\del-\del'=d-1.
\end{align}

First, suppose that $Y$ does not contain $L$. Let $\card(Y\cap (L \setminus \Omega))=e$. Here $e\leq |\beta'|$.  Then by Lemma \ref{dimLemma}, 
\begin{align}
2(d-1)+g(Y)+e-1\geq  \dim V^{d,\delta}(\al,\bet)-2=(2d+g+|\bet|-1)-2.
\end{align}
Equivalently, $(g(Y)-(g-|\bet'-\bet|+1))+(e-|\bet'|)\geq -1$. On the other hand, neither of the summands on the left hand side is greater than $0$. If $g(Y)=g-|\beta'-\beta|+1$ and $e=|\beta'|$, then $C$ is not a point of indeterminacy. Indeed, in this case, $C$ cannot be an image of a stable map of genus $g$ that contracts a component of a positive genus. 

Suppose $g(Y)= g-|\beta'-\beta|$. Then $e=|\beta'|$ which forces all inequalities in \eqref{dim-ineq} to be equalities and $Y$ to be unibranch at every point of $L \setminus \Omega$. Therefore, $Y$ is a generic point of $V^{d-1,\delta'+1}(\al',\bet')$, a generalized Severi variety satisfying condition \eqref{indeterminacy condition}.

Consider now the case when $e=|\beta'|-1$ and $g(Y)=g-|\beta'-\beta|+1$. If $Y$ is unibranch at every point of $(L \setminus \Omega)$, then by Lemma \ref{dimLemma}, the curve $Y$ belongs to the Severi variety $V^{d-1,\del'}(\alpha', \beta'')$, with $|\beta''|=|\bet'|-1$. Informally, we see two of the new points points of tangency on $V^{d,\del}(\al,\bet)$ coalesce. In a such situation, the newly formed tacnode of order $m$ on $Y\cup L$ is a limit of only $m-2$ nodes in the nearby fibers of $S^{d,\del}(\alpha,\beta)$. 

Suppose $Y$ is not unibranch along $L \setminus \Omega$. We consider a nodal reduction of a generic one-parameter family in $S^{d,\delta}(\al,\bet)$ with the center at $[C]$.
We use notations of Section \ref{stable reduction} throughout. Since $[C]$ contains $L$ with multiplicity $1$, we have $c_1=1$. Let $\beta_1$ is the number of sections in $V$ meeting $\overline{Y}$. We set $\beta_2=|\bet|-\beta_1$. Note that
$$\beta_1+\card((P'+P''_v)\cap \overline{Y})=\card(\pi^{-1}(L \setminus \Omega))\geq e+1 =|\bet'|.$$ 
Also, $\card( P''_f\cap \overline{Y})\geq f$ and $v\leq \beta_2$.
 Putting everything together we have the following inequalities
\begin{align*}
g &=p_a(Z_0)\geq g(\overline{Y})+(1-c_1-c_2)+\card(P'\cap \overline{Y})+\card(P'' \cap \overline{Y})-1 \\
&=g(\overline{Y})-(v+f)+\card(P'\cap \overline{Y})+\card(P'' \cap \overline{Y})-1 \\
 &= g+(\card(P''_f\cap \overline{Y})-f)+(\beta_1+\card(P'\cap \overline{Y})+\card(P_v''\cap \overline{Y})-|\beta'|)+(\beta_2-v) \geq g.
\end{align*}
We conclude that:
\begin{enumerate}
 \item\label{gen} All connected components of $P$ have arithmetic genus $0$. 
 \item\label{f} $\card(P''_f\cap \overline{Y})=f$.
 \item \label{card} $\card(\pi^{-1}(L \setminus \Omega))= e+1 =|\bet'|$.
 \item Connected component of $P''_v$ are in one to one correspondence with $\beta_2$ sections of $V$ not meeting $\overline{Y}$.
\end{enumerate}

\begin{figure}[htbp]
\begin{center}
\epsfig{file=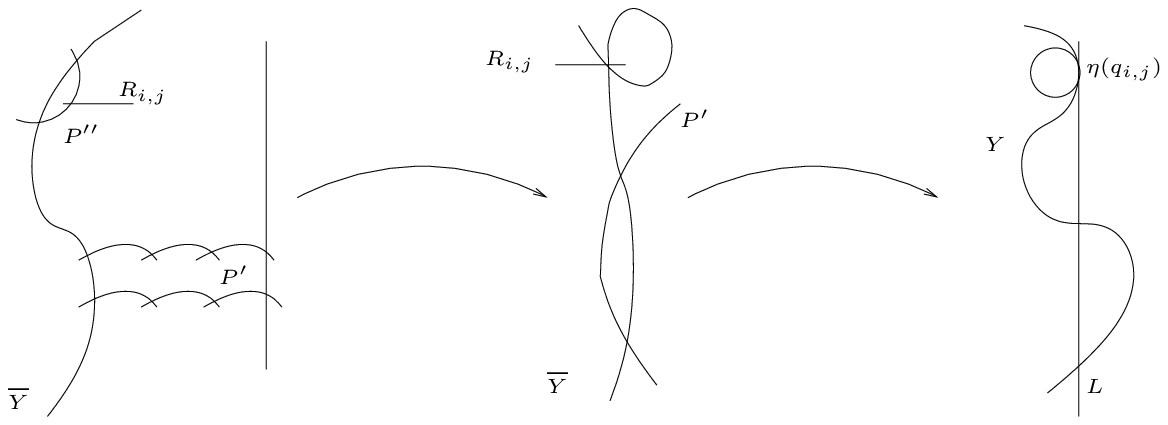, scale=0.8}
\caption{}
\label{fig1}
\end{center}
\end{figure}

\begin{figure}[htbp]
\begin{center}
\epsfig{file=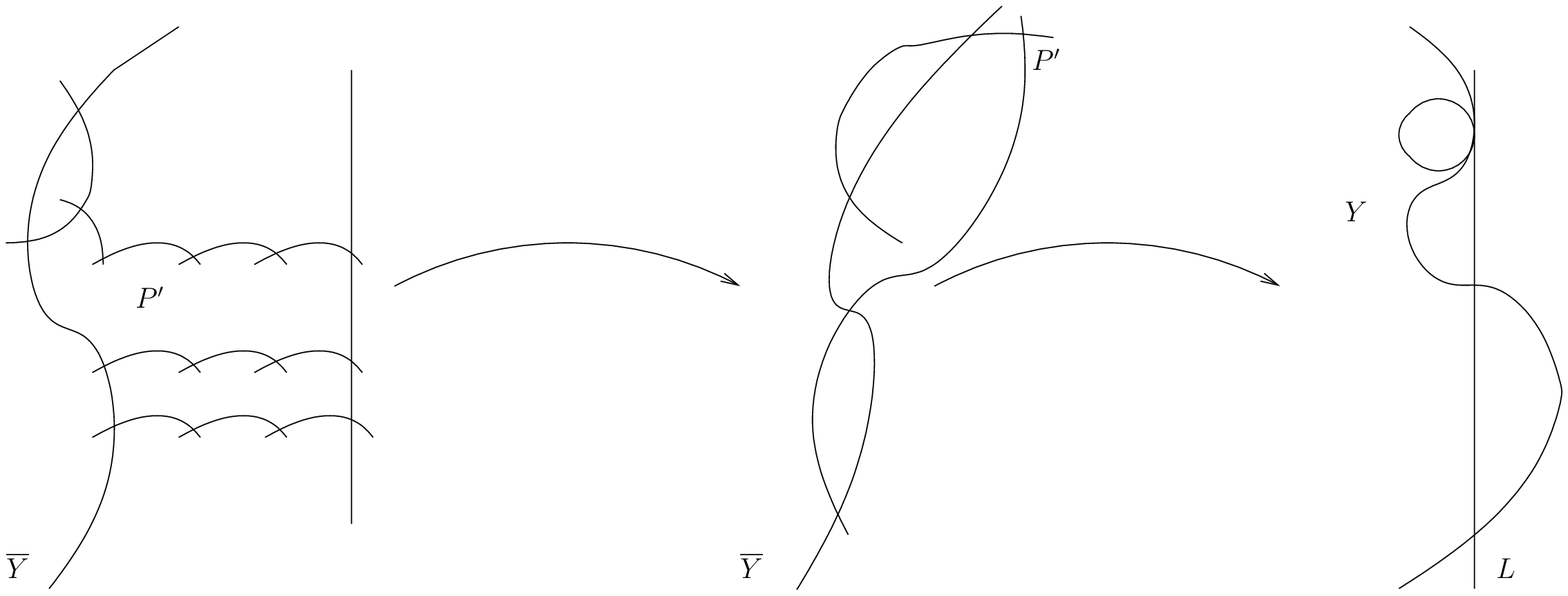, scale=0.4}
\caption{}
\label{fig2}
\end{center}
\end{figure}

It follows from \eqref{gen} and \eqref{f} that $P''_f$ is empty. Also, from \eqref{gen} any component of $P''_v$ has to meet $\overline{Y}$ in at least two points, and from \eqref{card} it follows that $v$ is at most $1$. In every case, the stable limit is one of the finitely many curves. We can have either $v=1$, with the only component of $P''_v$ being a $\pp^1$ meeting $\overline{Y}$ in two points and meeting one of the sections $R_{i,j}$. In this case, $P'$ meets $\overline{Y}$ in points which all map to different points on $L$. The other possibility is that $P''$ is empty and $P'$ meets $\overline{Y}$ in a set of points among which there are two that map to the same point on $L$. The first possibility is presented in Figure \ref{fig1}, where the curve on the left is a central fiber of the nodal reduction, the curve in the center is its stabilization. The curve on the right is the planar image. The second possibility is presented in Figure \ref{fig2}. It follows that $[C]$ is not a point of indeterminacy.


{\bf Remaining case: $Y$ contains $L$.}
We set $Y=X\cup L$ where $Y$ is a general point of the Severi variety $V^{d,\del''}(\al'',\beta'')$ with $|\beta''-\bet'|+\del-\del''=d-2$. Together with \eqref{CND}, this is equivalent to $$g(X)=g-|\bet''-\bet|+2.$$

Following the notations of the Section \ref{stable reduction}, we consider the nodal reduction of a generic family in $S^{d,\delta}(\alpha,\beta)$ with the center at $[C]$. 
 
Note that $c_1\leq 2$ by degree considerations, and $(P''\cap \overline{X}) \geq c_2$ since $Z_0$ is connected. We have
\begin{align*}
g &=p_a(Z_0)\geq g(\overline{X})+(1-c_1-c_2)+(P'\cap \overline{X})+(P''\cap \overline{X})-1 \\
 &\geq g-|\bet''-\bet|+(P'\cap \overline{X})+(2-c_1)+((P''\cap \overline{Y})-c_2)\geq g.
\end{align*}

This is possible only if all the inequalities are equalities. We draw the following conclusions.
\begin{enumerate}
\item All connected components of $P$ have arithmetic genus $0$. 
\item $P''$ is empty. 
\item There are two components, $P'_1$ and $P'_2$, of $P'$. Each is a tree of rational curves that is mapped with degree $1$ onto $L$. 
\end{enumerate}

\begin{figure}[htbp]
\begin{center}
\epsfig{file=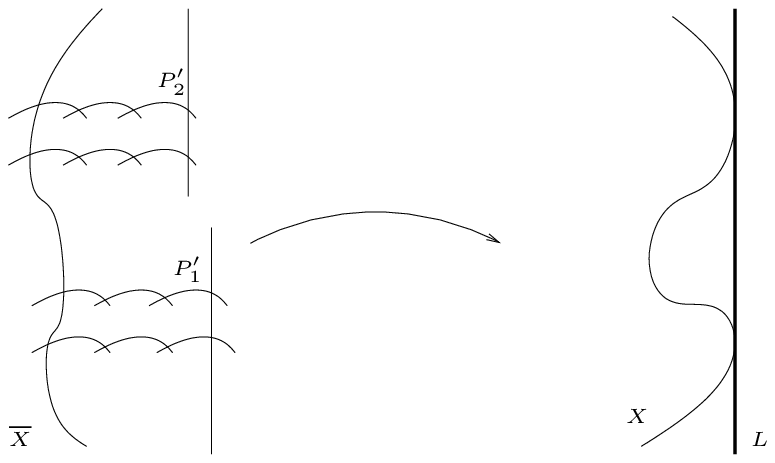, scale=0.6}
\caption{}
\label{fig3}
\end{center}
\end{figure}

We conclude that the stable reduction of $Z_0$ looks like the curve in Figure \ref{fig3}. In particular, there are finitely many possible stable limits of one-parameter families with the center at $[C]$, and so
$[C]$ is not a point of indeterminacy. 

The proof of Proposition \ref{points of indeterminacy} is finished.

We conclude this chapter with the restatement of Equation \eqref{L-induction}:
\begin{equation}\label{L-induction2}
\begin{split} 
L^{d,\delta}(\alpha,\beta) &= \sum_k kL^{d,\delta}(\alpha+e_k,\beta-e_k) \\
&+\sum I^{\beta'-\beta}\binom{\alpha}{\alpha'}\binom{\beta'}{\beta}L^{d-1,\delta'}(\alpha',\beta') 
+\discr{S^{d,\delta}(\alpha,\beta)}{H_q}{\lambda}.
\end{split}
\end{equation}

\section{Local Geometry of $\Y^{d,\delta}(\al,\beta)$}\label{local geometry}
In this chapter, we describe the geometry of the family $\Y^{d,\delta}(\al,\bet)\ra S^{d,\delta}(\al,\bet)$. Note that the only singularities, besides nodes, that appear on the curves corresponding to generic points of Type II components, and on the curves corresponding to the points of indeterminacy, are higher-order tacnodes. Therefore the description of the local geometry of $\Y^{d,\delta}(\al,\beta)$ invariably invokes versal deformation spaces of (arbitrary order) tacnodes. We recall the necessary definitions in what follows.

\subsection{Versal deformation space of a tacnode}\label{deformation space}

Recall that the versal deformation space of the $m^{\th}$-order tacnode $y^2-x^{2m}=0$ is $T\cong \cc^{2m-1}$. We let the coordinates on $T$ be $(a_2,\dots, a_{2m})$. Then the miniversal family $\Y \ra T$ can be defined by the equation
\begin{align}
\label{equation-of-versal}
y^2=(x^m+a_2 x^{m-2}+a_3 x^{m-3}+\cdots+a_m)^2+a_{m+1}x^{m-1}+\cdots+a_{2m-1}x+a_{2m}
\end{align}
inside $T\times \spec \cc[x,y]$.

For $a=(a_2,\dots, a_{2m})\in T$, we write $\Psi(a_2,\dots, a_{2m})(x)$, or $\Psi_a(x)$, to denote the polynomial on the right hand side of Equation \eqref{equation-of-versal}. We use $\Psi_a(x,z)$ to denote the homogenization of $\Psi_a(x)$. 

Let $D\cong \cc^{2m-1}$ be the space of monic polynomials of degree $2m$ with a trivial $x^{2m-1}$ coefficient. Any polynomial in $D$ can represented as a sum of the square of a monic polynomial of degree $m$ and a polynomial of degree $m-1$. In fact, 
the map $\Psi: T \ra D$ 
sending a point $a\in T$ to a polynomial $\Psi_a(x)$ is an isomorphism, since its Jacobian is an upper-triangular matrix with non-zero complex numbers along the diagonal. 

\begin{definition} A {\bf combinatorial type} of a polynomial $\Psi$ of degree $d$ is an $r$-tuple $(m_1,\dots,m_r)$ of multiplicities of distinct roots of $\Psi$. Here $r$ is the number of distinct roots of $\Psi$.
\end{definition}

The deformation space $T$ has a natural geometric stratification given by the combinatorial type of the polynomial $\Psi_a(x)$.

\begin{definition} We denote 
the stratum of deformations of type $(m_1,\dots,m_r)$ by $$\Delta^0\{m_1,\dots,m_r\};$$
and denote its closure by $\Delta\{m_1,\dots,m_r\}$. 
\end{definition}

\begin{remark}
For $a \in \Delta^0\{m_1,\dots,m_r\}$, the singularities of the fiber $y^2=\Psi_a(x)$ are double points $y^2=x^{m_i}$, therefore $\Delta^0\{m_1,\dots,m_r\}$ is an equisingular stratum.  
 \end{remark}
We also define $$\Delta_r:=\Delta\{\underbrace{2,2,\dots,2}_{r},\underbrace{1,\dots,1}_{2m-2r}\}$$
to be the closure of the locus of $r$-nodal curves.

 We will distinguish the hyperplane $H_{2m}$ inside $T$ defined by $a_{2m}=0$. Under the identification of the tangent space 
$\mathbb{T}_0 T$ with the space $\Def_1(y^2=x^{2m})$ of the first-order deformations of the tacnode,
$\mathbb{T}_0 H$ corresponds to the first order deformations of $y^2=x^{2m}$ vanishing at $(0,0) \in \cc^2$. 

\subsection{Multiple tacnodes}
 Let $\m=(m_1,\dots,m_n)$ be a sequence of positive integer numbers. For $1\leq i\leq n$, let $T(i) \cong \cc^{2m_i-1}$ be the versal deformation space of the tacnode $y^2-x^{2m_i}=0$, and let $(a_{i, 2},\dots,a_{i, 2m})$ be coordinates on $T(i)$. 
We set
\begin{align*}
T:=\prod_{i=1}^{n} T(i) \cong \cc^{2(\sum_i m_i)-n}
\end{align*}
to be the product of these deformation spaces. 

By analogy with a single tacnode case, we define the following loci inside T:
\begin{align*}
\Delta_{\m} &:=\prod_{i=1}^n \Delta_{i,m_i} \, , \\
\Delta_{\m-\mathfrak{1}} &:=\prod_{i=1}^n \Delta_{i,m_i-1} \, , \\
\Delta_{\m-\mathfrak{1}-e_j} &:= \prod_{i=1}^{n} \Delta_{i,m_i-1-\delta_{i,j}}, \quad  1\leq j \leq n.
\end{align*}

\subsection{Geometry at the generic point of Type II component}\label{generic point}
We recall the description, given in \cite[ Section 4.4]{CH1}, of the geometry of $S^{d,\delta}(\al,\bet)$ around a generic point of a Type II component of $H_q$.

Suppose that $[X]=[Y\cup L]\in  C_L^{d-1,\delta'}(\alpha',\beta')$ is a generic point. Here, $[Y]$ is a generic point of $C^{d-1,\delta'}(\alpha',\beta')$. The curve $X$ has $|\alpha'|+|\beta'|$ tacnodes, at the points of contact of $Y$ with $L$, of which $|\alpha'|$ are in the set $\Omega \subset L$. 
Among the ``moving'' $|\beta'|$ points of contact, exactly $|\beta|$ are the limits of  ``moving'' points of tangency in the nearby fibers. We call the corresponding tacnodes of $X$ ``old''\footnote{This does not mean that an old tacnode is a limit of tacnodes in the nearby fibers.}. 

The rest of the tacnodes correspond to ``new'' points of tangency of $Y$ with $L$.
Set $\beta'':=\beta'-\beta$ and $n:=|\beta''|$. Denote the new tacnodes of $X$ by $y_1,\dots, y_n$, and their multiplicities 
by $m_1,\dots, m_n$. Consider the Severi variety $V^{d,\delta''}(\al,\bet)$, where
$\delta''=\delta-(I\beta''-n)$. Inside it, we define an open subvariety $V$ consisting of the deformations of $X=Y\cup L$ in $V^{d,\delta''}(\al,\bet)$ satisfying the following conditions:
\begin{enumerate}
\item Tangencies at $\Omega$ are preserved.
\item The deformations of the $|\beta|$ ``old'' tacnodes preserve two branches: the deformation of $Y$ and $L$, respectively.
\item The deformation of $Y$ near an ``old'' tacnode of order $i$ remains tangent (at an unspecified point) to the line $L$ with the multiplicity $i$. 
\end{enumerate}
The subvariety $V$ is a {\it relaxed Severi variety} in the sense of \cite[Proposition 4.8]{CH1}. Then, in the neighborhood of $X$, the variety $V^{d,\delta}(\al,\bet)$ is a closure of deformations of $X$ inside $V$ such that every tacnode of order $y_i$ deforms to $m_i-1$ nodes. 

To say it differently, consider a map
$$\phi: V \ra T:=\prod_{i=1}^n \Def(y_i, X).$$
Let $\{a_{i,j}\}_{2\leq j \leq 2m_i}$ be the coordinates on $T_i:=\Def(y_i,X)$ as described in Section \ref{deformation space}.
Then summarizing the results of \cite[Section 4.4]{CH1}, and Lemma 4.9 \cite{CH1} in particular, we have: 
\begin{enumerate}
\item $V$ is smooth.
\item $\phi$ is smooth at $[X]$.
\item $W:=\phi(V)$ contains $\Delta_{\m}$ and is smooth of dimension $\sum(m_i-1)+1$.
\item The tangent space to $W$ at the origin is not contained in the union of hyperplanes $a_{i,2m_i}=0$.
\item $W\cap \Delta_{\m -\mathfrak{1}}=\Delta_{\m} \cup \Gamma_{\m}$, where $\Gamma_{\m}$ is a curve intersecting $\Delta_{\m}$ at the origin with the multiplicity 
$$M=\prod_{i=1}^n m_i.$$
\item $V^{d,\delta}(\al,\beta)=\phi^{-1}(\Gamma_{\m})$ and $V_L^{d-1,\delta'}(\al',\bet')=\phi^{-1}(0)$.
\end{enumerate}

We note in passing, the above discussion implies that the stable reduction of the family $$\Y\ra S^{d,\delta}(\al,\beta)$$ at a generic point $[X]=[Y\cup L]$ of a Type II component $C_L^{d-1,\delta'}(\al',\beta')$ involves no blow-ups. Indeed, since the ``new" tacnode $y_i$, of order $m_i$, is a limit of $m_i-1$ nodes in the nearby fiber, the stable limit of any one-parameter family with the center at $[X]$ will have a node lying over $y_i$. All other singularities are resolved in the process of the stable reduction. We conclude that the stable limit will be a nodal union of $Y^{\nu}$ with the line $L$ at the points lying over $y_i$'s. 

\subsection{Geometry at the point of indeterminacy}\label{geometry indeterminacy}

We retain the notations of the previous section.
Suppose $[X]=[Y\cup L]$ is a point of indeterminacy, where $Y$ is a generic point of $V^{d-1,\delta'}(\al',\bet')$ satisfying $|\bet'-\bet|+\delta-\delta'=d-2$. Note that there are precisely 
\begin{align}\label{number-of}
\binom{\alpha}{\alpha'}N^{d-1,\delta'}(\alpha',\beta')
\end{align} points of indeterminacy corresponding to the triple $(\delta',\al',\beta')$.

First, we note that $S^{d,\delta}(\alpha,\beta)$ has several analytic branches in the neighborhood of $[X]$. Each branch is specified by the choice of tacnodes of $[X]$ which are the limits of the $|\beta|$ moving points of tangency in the nearby fibers. There are $\binom{\beta'}{\beta}$ such branches. Choose one of them, and denote it by $S$. Let $y_1, \dots y_n$ be the remaining ``new" tacnodes of $X$, of order $m_1, \dots, m_n$, respectively (here, $n=|\beta'-\beta|$). Let $V$ be the {\it relaxed Severi variety} of the previous section, but now defined inside the Severi variety $V^{d,\delta''}(\al,\bet)$, where $\delta''=\delta-(I\beta''-n)+1$. Then, by \cite[Lemma 4.9]{CH1}, the map
$$\phi: V \ra T:=\prod_{i=1}^n \Def(y_i, X),$$
satisfies the following conditions:
\begin{enumerate}
\item $V$ is smooth.
\item $\phi$ is smooth at $[X]$.
\item $W:=\phi(V)$ contains $\Delta_{\m}$ and is smooth of dimension $\sum(m_i-1)+1$.
\item The tangent space to $W$ at the origin is not contained in the union of hyperplanes $a_{i,2m_i}=0$.
\end{enumerate} 
However, the geometry of $V^{d,\delta}(\al,\bet)$ in the neighborhood of $[X]$ in $V$ is more complicated. We proceed now to describe it. 

Speaking colloquially, we see only $\sum_{i=1}^{n}(m_i -1)-1$ nodes approaching tacnodes $y_1,\dots, y_n$. Hence, there is a tacnode $y_{i_0}$ that is the limit of only $m_{i_0}-2$ nearby nodes. Every other tacnode $y_i$ ($i\neq i_0$) is the limit of $m_i-1$ nearby nodes.  

To make the above statement precise, we consider, for each $1\leq i \leq n$, the intersection $$W\cap \Delta_{\m-\mathfrak{1}-e_i} =\Delta_{\m} \cup S_{i}, $$
where $S_i:=\overline{(W\cap \Delta_{\m-\mathfrak{1}-e_i})\setminus \Delta_{\m}}$ is residual to $\Delta_{\m}$ in the intersection.
We set $$S_\m:=\bigcup_{i=1}^n S_{i}.$$
Then $S$ is identified with $S_{\m}$ via the map $\phi$. The hyperplane section $H_q$ of $S$ is identified with $S_{\m}\cap \Delta_{\m}$. Symbolically,
$$S=\phi^{-1}(S_{\m}) \ \ \text{and} \ \ S\cap H_q=\phi^{-1}(S_{\m}\cap \Delta_{\m}).$$ 
By the definition of $\phi$, locally around the point $y_i$ in the central fiber $X$, the family $\Y$ is isomorphic to the pullback of the miniversal family over $\Def(y_i,X)$. 

We thus reduce the calculation of the discrepancy along $H_q$ on the branch $S$ of $S^{d,\delta}(\alpha,\beta)$ to the calculation of the discrepancy along the divisor $\Delta_{\m}$ on the surface $S_\m$ inside the product of the versal deformation spaces of tacnodes. Recalling that the number of the points of indeterminacy of the type $(\delta',\al',\beta')$ is given by \eqref{number-of}, and that there are $\binom{\beta'}{\beta}$ choices for the branch $S$, we conclude that
\begin{align}\label{discrepancy-total}
\discr{S^{d,\delta}(\al,\beta)}{H_q}{\lambda}= \sum \binom{\al}{\al'}\binom{\bet'}{\bet}N^{d-1,\delta'}(\al',\bet')\cdot \discr{S_{\m}}{\Delta_\m}{\lambda},
\end{align}
where the sum is taken over the triples $(\delta',\al',\bet')$ satisfying $|\bet'-\bet|+\delta-\delta'=d-2$ and $$\m=\{\underbrace{1,\dots,1}_{\beta'_1-\beta_1}, \underbrace{2,\dots,2}_{\beta'_2-\beta_2}, \dots\}.$$

\section{Geometry of the versal deformation space of tacnode}\label{single tacnode}

In the previous chapter, we saw that the local geometry of the Severi variety is reflected in the geometry of certain loci in the product of the versal deformation spaces of tacnodes. 
In this chapter, we lay out the framework for studying the deformation space of a single tacnode. We obtain several results which elucidate the geometry of the miniversal family. They will be generalized to multiple tacnodes in the following chapter, and used to in the calculation of discrepancies. 

  
\subsection{Alterations}
We use the notations of Section \ref{deformation space}. Recall that the versal deformation space of the singularity $y^2=x^{2m}$ is $T=\spec \cc[a_2,\dots,a_{2m}]$. The miniversal family $\Y$ over $T$ is defined by
\begin{align}
\label{equation-of-versal-2}
y^2=\Psi_a(x)=(x^m+a_2 x^{m-2}+a_3 x^{m-3}+\cdots+a_m)^2+a_{m+1}x^{m-1}+\cdots+a_{2m-1}x+a_{2m} 
\end{align}
inside  $T\times \spec \cc[x,y]$.

We make a base change $\pi_1\co T'=\spec \cc[b_2, b_3,\dots, b_{2m}]\ra T$ defined by $a_i=b_i^i$. For 
$$\mu:=\mu_{2}\times \mu_3\times \dots \times \mu_{2m},$$ 
where $\mu_r$ is the cyclic group of $r^{\th}$ roots of unity, we have $T=T'/\!/ \mu$.

We let $\pi_2\co T''\ra T'$ to be the blow-up of $T'$ at the origin, whose exceptional divisor we denote by $E$, and set $$\pi:=\pi_2\circ\pi_1\co T''\ra T.$$
Denote by $\Y''\subset T''\times \spec \cc[x,y]$ the pullback of $\Y$ to $T''$. Let $\I_E$ be the ideal sheaf of $E$ on $T''$, and consider the ideal sheaf 
$$\I:=((\I_E, x)^m, y)$$ on $T''\times \spec \cc[x,y]$.

First, note that $Bl_{\I} (T''\times \spec \cc[x,y]) \ra T''$ is a family of surfaces with fibers over $T'' \setminus E$ being affine planes $\cc^2$, and fibers over the points in the exceptional divisor  $E$ being the union of $Bl_{(x^m,y)}\cc^2$ and 
$\pp(1,1,m)=\proj \cc[x,y,z]$.  
Here, $z$ stands for a local generator of $I_E$ and $\cc[x,y,z]$ is graded with $\deg x=\deg z=1$ and $\deg y=m$. Finally, we denote 
\begin{align*}
\Z:=Bl_{\I} \Y''   \  \ \text{and}  \ \ F \co \Z \ra T''.
\end{align*} 

The family $F\co \Z\ra T''$ can be seen as a first step towards the stable reduction of the miniversal family $\Y\ra T$. The following paragraphs make this statement more precise.

Consider a usual open cover of the blow-up $T''$ by the affine charts 
$$D_{(b_i)} T'' :=\spec \cc[b_2/b_i,b_3/b_i,\dots, b_{i},\dots, b_{2m}/b_i].$$
Then, over $D_{(b_i)}T''$, the family $\Z$  has equation
$$y^2=\Psi\biggl(\left(\frac{b_2}{b_i}\right)^2,\dots, \left(\frac{b_{2m}}{b_i}\right)^{2m}\biggr)(x,b_i).$$ 
The restriction $E_{(b_i)}$ of the exceptional divisor to $D_{(b_i)}T''$ is given by equation $b_i=0$, and the restriction of $\Z$ to $E_{(b_i)}$ is 
$$\Z_E :=\Z \times_{T''} E_{(b_i)}=\mathcal{S}\cup \T,$$
where $\mathcal{S}$ is defined by
$$\{(y/x^m)^2=1\}\subset E_{(b_i)}\times Bl_{(x^m,y)}\cc^2$$ and $\T$ is defined by
$$ \{y^2=\Psi\biggl(\left(\frac{b_2}{b_i}\right)^2,\dots, \left(\frac{b_{2m}}{b_i}\right)^{2m}\biggr)(x,z) \} \subset E_{(b_i)}\times \pp(1,1,m).$$

In words, over $E_{(b_i)}$, the family $\Z_E$ has two components. The first component, $\mathcal{S}$, is a trivial family with the normalization of $y^2=x^{2m}$ as its fibers. The second component, $\T$, is a family of divisors, not passing through the vertex, in the linear series $|\O(2)|$ on $\pp(1,1,m)$\footnote{The weighted projective space $\pp(1,1,m)$ is a projective cone over the rational normal curve of degree $m$. The line bundle $\O(1)$ is the restriction of the hyperplane section.}. We call fibers of $\T$ {\bf tails}. Note that all such tails are hyperelliptic curves of arithmetic genus $m$ with $2$ distinguished mark points defined by $z=0$. The marked points are the points of intersection with $\mathcal{S}$ and are exchanged by the hyperelliptic involution $y\mapsto -y$.

We observe that singularities appearing on the hyperelliptic tails are planar double points $y^2=x^n$ with $n<2m$, and so are strictly better singularities than the tacnode $y^2=x^{2m}$. Hence, after an alteration of the base and an appropriate blow-up of the total family, we arrived at the family of curves with milder singularities. Conceivably, by repeating the procedure for all singularities, we would arrive at the family of the stable curves. However, for our purposes, the family $F: \Z \ra T''$ will suffice.
 
{\it Caution:} Even though the family $\mathcal{S}$ is a trivial family when restricted to affine charts $E_{(b_i)}$, it is not trivial globally. 

\subsection{Local charts} 

Note that all blow-ups in the previous subsection were made with $\mu$-invariant centers. Hence, the action of $\mu$ extends to $\Z$ making the morphism $F$ equivariant. 
For every $2\leq i \leq 2m$, we consider a quotient of $F\co \Z \ra T''$ by the natural action 
of $$\tilde{\mu}_i := \mu_2\times \cdots \times \hat{\mu}_i \times \cdots \times \mu_{2m}.$$

We set $\Z_i:=\Z /\!/ \tilde{\mu}_i$ and $T_i:=T'' /\!/ \tilde{\mu}_i$. The quotient morphism $F_i \co \Z_i \ra T_i$ is nothing else than the weighted blow-up of $T$ with weights $(2,3,\dots,2m)$, followed by a base change $a_i=b_i^{i}$. We define $\T_i$ to be the quotient $\T /\!/ \tilde{\mu}_i$ of the tails component.
By abuse of notation, we also use $E$ to denote the exceptional divisor of $F_i$. If $W$ is a subvariety of $T$, we denote by $W_E$ its exceptional divisor in $T_i$.

Note that the action of $\tilde{\mu}_i$ is free on $D_{(b_i)}T''$. The quotient, denoted $D_{(b_i)}T_i$, is isomorphic to $\spec \cc[c_2,\dots, b_i,\dots, c_{2m}]$, via $c_j=a_j/b_i^{j}$. The equation of $E_{(b_i)}:=E\cap D_{(b_i)}T_i$ is $b_i=0$. 
Over $E_{(b_i)}$, the equation of $\T_i$  is
$$\{y^2=\Psi(c_2,\dots,c_{i-1},1,c_{i+1},\dots,c_{2m})(x,z) \} \subset E_{(b_i)}\times \pp(1,1,m).$$
Henceforth, in our discussion we will identify the family $\T_i \ra E_{(b_i)}$ with the affine space $$\spec[c_2,\dots,\hat{c}_i,\dots, c_{2m}]$$
 of in-homogeneous polynomials 
$$\Psi(c_2,\dots,c_{2m})(x)=x^{2m}+c_2 x^{2m-2}+\cdots+c_{2m}$$ of
degree $2m$ with $c_i=1$. 

\begin{lemma}\label{versality}
Given any point $p \in E_{(b_2)}$, denote by $\D(p)$ the product of the deformation spaces of the singularities of the fiber $(\T_2)_p$. Then the family $\T_2 \ra E_{(b_2)}$ induces a smooth morphism from an analytic neighborhood of $p$ in $E_{(b_2)}$ to an analytic neighborhood of the origin in $\D(p)$.
\end{lemma}

\begin{proof}
Suppose $p=(c_3,\dots,c_{2m}) \in \Delta^0\{m_1,\dots, m_r\}$. Set $P:=\Psi(1,c_3,\dots,c_{2m})=\prod_{i=1}^{r} (x-x_i)^{m_i}$, where $x_i$ satisfy $\sum_{i} x_i m_i=0$.  The equation of $(\T_2)_p$ is $y^2=P(x,z)$ and so the singularities of $(\T_2)_p$ are double points $y^2=x^{m_i}$. Hence, $\D(p)\cong \cc^{\sum (m_i-1)}$, with equisingular locus being $0 \in \D(p)$. 

The first order deformations of $p$ correspond to the first order deformation 
of $P$ of the form $P+\epsilon Q$, with $Q$ being an arbitrary polynomial of degree $2m-3$.
Equisingular deformations are precisely those satisfying 
$$\prod_{i=1}^{r}(x-x_i)^{m_i-1} \,  | \, Q(x).$$
Lemma follows from the fact that this divisibility condition imposes $\sum_{i=1}^{r} (m_i -1)$ independent linear conditions on the coefficients of $Q$.
\end{proof}

\subsection{Tangent cones} 
\begin{lemma}\label{tang_cone} Suppose $W$ is a subvariety of $T$ that is smooth at the origin, has dimension $m$, contains $\Delta_{m}$, and whose tangent space is not contained in the hyperplane $H_{2m}$. Then the exceptional divisor of $W$ in $T''$ is given by equations 
$$a_{m+i}=b_{m+i}^{m+i}=0, \ 1\leq i\leq m-1.$$
Here $[b_2:\ldots:b_{2m}]$ are homogeneous coordinates on $E\subset T''$.
\end{lemma}

\begin{proof}
The tangent space of $W$ in $T$ is a linear space of dimension $m$ not contained in the hyperplane $a_{2m}=0$, but containing $\{a_{m+1}=\dots=a_{2m-1}=a_{2m}=0\}$. Therefore, the initial ideal of $W$ satisfies
\begin{align*}
\textbf{in}(I(W))=\ker [\begin{array}{c|c|c} \mathbf{0}_{m-1} & \mathbb{I} & * \end{array}],
\end{align*}
where $\mathbf{0}_{m-1}$ is the zero $(m-1)\times (m-1)$ matrix, $\mathbb{I}$ is a row-permutation of the identity matrix $\mathbb{I}_{m-1}$, and $*$ is a column of complex numbers. Equivalently,
$$\textbf{in}(I(W))=(a_{m+\sigma(i)}+t_{i}a_{2m}), \ 1\leq i\leq m-1, \ t_i\in \cc, \ \sigma\in \mathfrak{S}_{m-1}.$$
Therefore, under the substitution $a_{i}=b_{i}^{i}$, the initial ideal of $\pi_1^{-1}(W)$ in $T'$ becomes 
$$\{b_{m+i}^{m+i}=0, \ 1\leq i\leq m-1\}.$$
\end{proof}

Starting with the simple Lemma \ref{tang_cone}, we draw important corollaries regarding various geometric strata inside $T$. Consider a point $a\in \Delta^0\{2m_1, \dots, 2m_n\}$ with $$\Psi_a(x)=(x-x_1)^{2m_1}\cdots (x-x_n)^{2m_n}$$ where $\sum_{i=1}^n m_i x_i=0$. We set $T(i)=\spec \cc [c_{i,j}]_{2\leq j \leq 2m_i}$
 to be the versal deformation space of the $m_i^{\th}$\nobreakdash-order tacnode $(x,y)=(0,x_i)$ in the fiber $\Y_a$. Then we have a natural map 
 $$\Phi: (T,a) \ra \prod_{i=1}^n T(i),$$
defined in a neighborhood of the point $a\in T$.

\begin{lemma}\label{versality-differential}
Consider the linear subspace  
$$W:=\{a_{m+1}=\dots=a_{2m-1}=0\}$$
of $T$ and a point $a\in W$. Then 
$d\Phi(\mathbb{T}_a W)$ is not contained in the union of hyperplanes $c_{i,2m_i}=0$. 
\end{lemma}
\begin{remark} We can also reformulate the corollary as the statement that, for all $i$, 
$$ (\mathbb{T}_a W)\cap (d\Phi)^{-1}(\{c_{i,2m_i}=0\})=\mathbb{T}_a \Delta_m.$$
\end{remark}
\begin{proof}
The statement follows from the geometric interpretation of the hyperplane $c_{i,2m_i}=0$ as the tangent space to deformations vanishing at the point $(x,y)=(x_i,0)$. By assumption, $\Psi_a(x)=P^2(x)$, for some polynomial $P(x)$. We note that 
the generic first order deformation $(P(x)+\epsilon Q(x))^2+\epsilon \lambda$ of $\Psi_a$ in $W$ does not vanish at $(x_i,0)$ for all $i$. The proof is finished. 
\end{proof}

\begin{corollary}\label{cor-W} Consider $W_E:=\{b_{m+1}=\dots=b_{2m-1}=0\} \subset E$. The following statements hold 
\begin{enumerate}
\item $W_E \cap (\Delta_{m-1})_E$ is analytically irreducible at every point of $(\Delta_m)_E$.
\item $W_E \cap (\Delta_{m-2})_E$ is analytically irreducible at all points of strata $\Delta^0\{2m_1,\dots, 2m_n\}_E$ where $(m_1,\dots, m_n)$ is an arbitrary $n$-tuple with $n\geq 3$.
\end{enumerate}
\end{corollary}

\begin{proof}

We first observe that both statements are equivalent to each of the analogous statements for $T''\setminus E$, $T' \setminus \{0\}$ and $T \setminus \{0\}$, in turn. 

Working now on $T\setminus \{0\}$, we consider a point $a\in  \Delta\{2m_1,\dots, 2m_n\}$ and the induced 
map
 $$\Phi: (T,a) \ra \prod_{i=1}^{n} T(i).$$

We observe that
 $$\Delta_{m-1}=\bigcup_{i=1}^n \Phi^{-1}(\prod_{j=1}^n \Delta_{m_j-\delta_{i,j}}).$$ 
 By  Lemma \ref{versality-differential}, we have $\mathbb{T}_a \Phi^{-1}(\prod_{j=1}^n \Delta_{m_j-\delta_{i,j}}) \cap \mathbb{T}_a W =\mathbb{T}_a \Delta_m$. Therefore, $\Delta_m$ is, locally at $a$, a smooth component of the intersection $\Phi^{-1}(\prod_{j=1}^n \Delta_{m_j-\delta_{i,j}}) \cap W$. This establishes the first claim. 
 
 The second claim is proved analogously. We observe that for $n\geq 3$, the locus $\Delta_{m-2}$, locally at $a$, is a union of varieties, each mapping onto $\Delta_{m_i}$ under $\Phi$ followed by the projection to some $i^{\th}$ factor. 
\end{proof}

\subsection{Geometry of $\Delta_{m-1}$}

We recall and reformulate  Lemma 2.12 of \cite{CH2} in the following form:
\begin{lemma}\label{poly}
Suppose $\lambda$ is a non-zero number.
For every positive integer $m$, there is a polynomial $P_m(x)=x^m+\alpha_{2}x^{m-2}+\dots+\alpha_m$ such that
$P_m(x)^2-\lambda$ has $m-1$ double roots. Moreover, $\alpha_{2k-1}=0$ for all $k$ , $\alpha_{2}\neq 0$, and $P_m(x)$ is unique up to scaling $\alpha_{2k}\mapsto \xi^k \alpha_{2k}$, where 
$\xi$ is an $m^{\th}$ root of unity.
 \end{lemma}

We now reprove Lemma 4.1 of \cite{CH1}.
\begin{lemma}\label{sectionG}
For any $m$-dimensional smooth variety $W$, containing $\Delta_m$, whose tangent space is not contained in the hyperplane $a_{2m}=0$, we have
$$W\cap \Delta_{m-1}=\Delta_{m}\cup \Gamma$$ 
where $\Gamma$ is a smooth curve tangent to $\Delta_m$ with order $m$ at the origin.
\end{lemma}

\begin{proof}
Recall that the morphism $f: T_2 \ra T$ is the weighted blow-up of $T$ followed by the base change of order $2$. Consider the weighted projective tangent cones of $W$ and $\Delta_{m-1}$, denoted respectively $W_E$ and $(\Delta_{m-1})_E$, inside the exceptional divisor $E$. Lemma \ref{poly} implies that away from the locus $b_2=0$, the intersection of $W_E$ and $(\Delta_{m-1})_E$ is a single point $G$, at least set-theoretically.

From now on, we work on $E_{(b_2)}=\spec \cc[c_3,\dots,c_{2m}]$. By Lemma \ref{tang_cone}, 
$$W_E=\{ c_{m+1}=\cdots=c_{2m-1}=0\}$$ 
and hence $$\T_2 \times_{E} W_E= \{y^2=(x^m+x^{m-2}+c_3 x^{m-3}+\cdots+c_m)^2+c_{2m}) \}.$$

Suppose coordinates of $G$ in $E_{(b_2)}$ are $(\lam_3,\dots, \lam_{2m})$, where $\lam_{2m}\neq 0$.
We then have
$$\Psi_{G}=P^2(x)+\lam_{2m}=Q^2(x)S(x)$$ 
where $P(x)=x^m+\frac{1}{2}x^{m-2}+\sum_{i=3}^{m}\lam_ix^{m-i}$,
the polynomial $Q(x)$ is monic of degree $m-1$, with distinct roots, and $S(x)$ is a quadric. We will proceed now to show that intersection $W_E \cap (\Delta_{m-1})_E$ is transverse at $G$. 

First, observe that the tangent space to $\Delta_{m-1}$ at $G$ is given by polynomials of degree $2m-3$ divisible by $Q$. The tangent space to $W_E$ at $G$ consists of the first-order deformations 
$$(P(x)+\epsilon R(x))^2+\lam_{2m}+\epsilon \lam'_{2m}=\Psi_G+\epsilon(2P(x)R(x)+\lam'_{2m}),$$
where $R(x)$ is a polynomial of degree $m-3$.

 To prove that the two tangent spaces intersect transversely, we need to show that $2P(x)R(x)+\lam'_{2m}$ is divisible by $Q$ only if $R=0$ and $\lam'_{2m}=0$. This is straightforward. Suppose $2P(x)R(x)+\lam'_{2m}$ is divisible by $Q$. Then
 \begin{align*}
 4P^2R^2+4\lam'_{2m}PR+(\lam'_{2m})^2 &\equiv 0  \bmod Q^2 \\
 -4\lam_{2m}R^2+4\lam'_{2m}PR+(\lam'_{2m})^2 &\equiv 0 \bmod  Q^2
 \end{align*}
 
 Observing that the left-hand side of the equality above has degree less that $2m-2=\deg Q^2(x)$,
 we conclude that 
 $$-4\lam_{2m}R^2+4\lam'_{2m}PR+(\lam'_{2m})^2=0.$$
 This implies that $\lam'_{2m}=0$ and $R=0$. 
 
We proved that $G=\overline{\Gamma}\cap E$ is a smooth point of the divisor $E \subset T_2$, Therefore $\overline{\Gamma}$ is a smooth curve, which in turn intersects the exceptional divisor $E$ transversely at $G$.

Let $H_2:=\{a_2=0\}$ and $H_{2m}:=\{a_{2m}=0\}$ be the coordinate hyperplanes in $T$. Then 
$$f^{*}H_2\cdot \overline{\Gamma}=(\overline{H}+2E)\cdot  \overline{\Gamma}=2E\cdot  \overline{\Gamma}=2.$$ Using the projection formula and the fact that $\overline{\Gamma}$ is a double cover of $\Gamma$, we deduce that $H_2\cdot\Gamma=1$. This necessarily implies that $\Gamma$ is smooth near the origin in $T$. Similarly, equalities 
$$f^* H_{2m}\cdot \overline{\Gamma}=(\overline{H}+2mE)\cdot\ \overline{\Gamma}=2mE\cdot \overline{\Gamma}=2m$$ imply that $H_{2m}\cdot \Gamma=m$. We finish by observing that $\Delta_m=W\cap H_{2m}$.

\end{proof}

It is enlightening to think of the family $\Z_{\, \overline{\Gamma}}:=\Z_{2} \times_{T_2} \overline{\Gamma}$ as a stable reduction of the family $\Y_{\Gamma}:= \Y \times_T \Gamma$. The central fiber of $\Z_{\, \overline{\Gamma}}$  
is the union of the normalization of $y^2=x^{2m}$ and the $(m-1)$\nobreakdash-nodal hyperelliptic tail $(\T_{2})_G$, while the generic fiber is an $(m-1)$\nobreakdash-nodal curve, by construction.
\begin{figure}[htbp]
\begin{center}
\epsfig{file=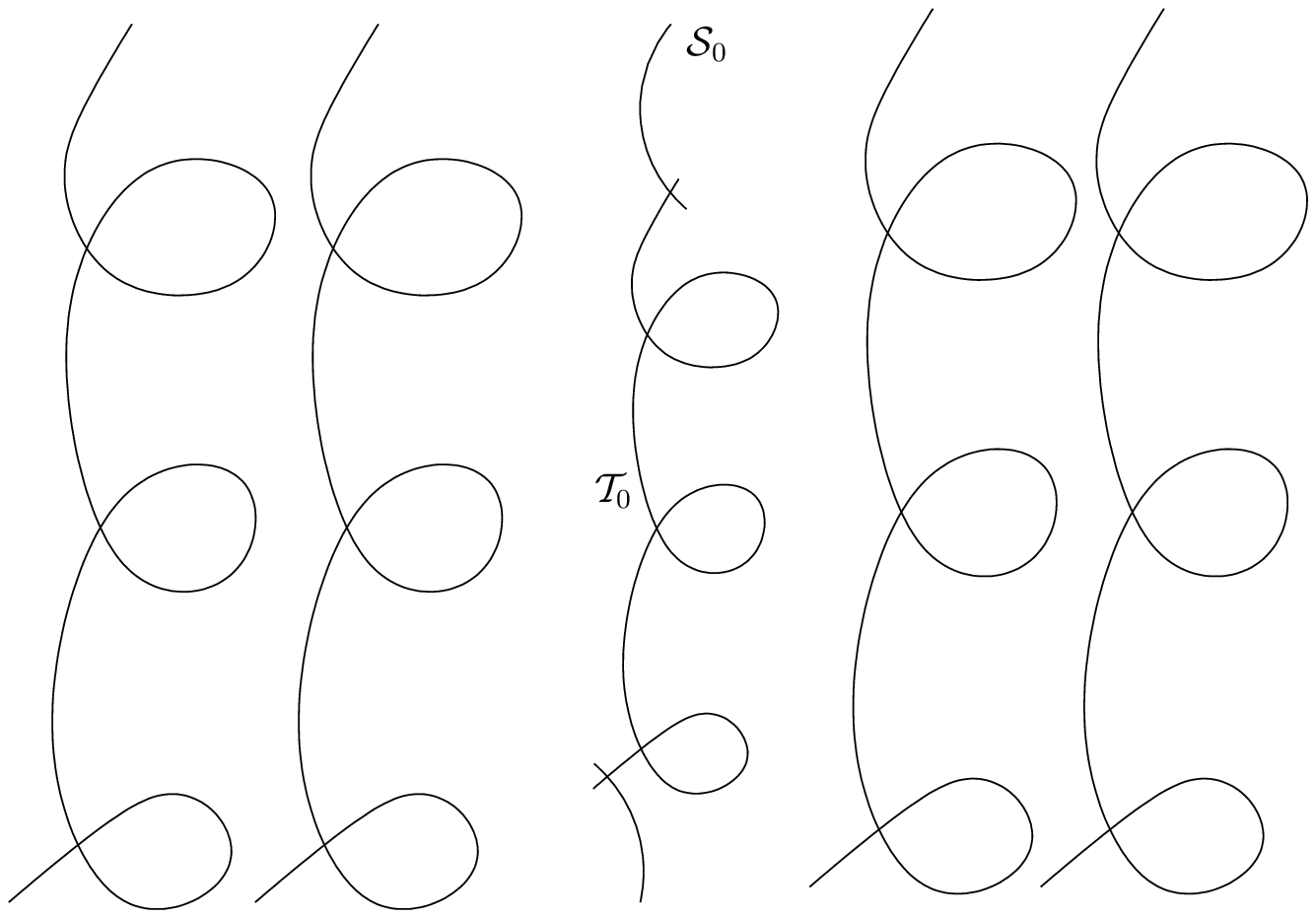, scale=0.45}
\caption{}
\label{fig4}
\end{center}
\end{figure}
It follows that the normalization of $\Z_{\, \overline{\Gamma}}$ is a family of generically smooth curves with the central fiber being the union of the normalization of $y^2=x^{2m}$ and the rational curve attached at the points lying over the tacnode of $\Y_0$. From the construction, we can easily deduce that the total space of the normalized family is smooth. It follows that any family of stable curves of genus $g$ that globalizes $(\Z_{\, \overline{\Gamma}})^{\nu}$ intersects the boundary $\Delta$ in $\Mg{g}$ with 
multiplicity $2$ at the point $G\in \overline{\Gamma}$. As $\overline{\Gamma}$ is a degree $2$ cover of $\Gamma$, ramified at $G$, we conclude that
any family of stable curves of genus $g$ that globalizes $(\Y_{\Gamma})^{\nu}$
intersects the boundary of $\Mg{g}$ with multiplicity $1$ at the origin. We state this more precisely as  
\begin{lemma}\label{delta-multiplicity}
The central fiber of the family $(\Y_{\Gamma})^{\nu}\ra \Gamma$ has a single node (lying over the tacnode in the central fiber of $\Y_{\Gamma}$). Moreover, the total space of the family is smooth at this node. 
\end{lemma}

\subsection{Geometry of $\Delta_{m-2}$}

Let $W$ be as above. Then we have $W \cap \Delta_{m-2}=\Delta_m \cup S$,
where $S$ is residual to $\Delta_m$ in the intersection. We would like to understand the geometry of $S$ and the restriction of the miniversal family to $S$. We do this by passing to $T''$. 

If $S_E$ is the exceptional divisor of $S$, then we have 
 $$W_E \cap (\Delta_{m-2})_E=(\Delta_m)_E \cup S_E.$$
 We study the geometry of the family $\T_S:=\T\times_E S_E$ of tails over $S_E$. First of all, we have the following result that follows from Corollary \ref{cor-W}:

 \begin{corollary}\label{DeltamS}
 \begin{align}
 S_E\cap (\Delta_m)_E=\bigcup_{i=1}^{\lfloor m/2\rfloor } \Delta\{2i,2(m-i)\}_E.
 \end{align}
 \end{corollary}

From the above discussion, we know that a generic fiber of $\T_S\ra S_E$ is an $(m-2)$-nodal hyperelliptic curve of arithmetic genus $m-1$ in the linear series $|\O(2)|$ on $\pp(1,1,m)$ and passing through points $[x:y:z]=[1:0:\pm 1]$. Hence, the normalization of a generic fiber is an elliptic curve, with two distinguished points. It follows
that there is an induced moduli map from the normalization $S_E^{\nu}$ to the coarse moduli space $\Mg{1,2}$. The map is defined, generically, by sending a point $p \in (S_E)^{\nu}$ to the normalization of the tail $(\T_S)_p$, marked at the points $[1:0:\pm 1]$ of attachment to $\mathcal{S}$. The first order of business is to understand which fibers of $\T_S$ are degenerate, i.e., have geometric genus $0$ or less.
 
\begin{proposition}\label{degenerate fibers}
The only degenerate fibers of $\T_S \ra S_E$ are fibers over the points 
\begin{align*} 
\Delta\{2i,2(m-i)\}_E \ \ \text{and} \ \ (\Delta_{m-1})_E\cap S_E.
\end{align*}
\end{proposition}

\begin{proof}
The $\delta$-invariant of the singularity $y^2=x^a$ is $\lfloor a/2\rfloor$. Therefore the geometric genus of a curve in a stratum $\Delta^0\{a_1,\dots,a_k\}_E$ is $m-1-\sum_i\lfloor a_i/2\rfloor$.

By Corollary \ref{DeltamS}, we have $S_E\cap (\Delta_m)_E=\bigcup \Delta\{2i,2(m-i)\}_E$. On the complement of $(\Delta_m)_E$, we have $b_{2m}\neq 0$. We next work on an open cover $U:=D_{(2m)}T_{2m}\subset T_{2m}$. 

Let $S_U=S_E\cap U$. Then any fiber over $S_U$ has equation $P(x)^2+1$, where 
\begin{align*}
P(x)&=(x^m+c_2 x^{m-2}+\cdots+c_m),  \ \ \text{(here $c_i=a_i/b_{2m}^i$}).
\end{align*} 
Suppose $P(x)^2+1$ has type $(m_1,\dots, m_k)$:
\begin{align*}
P(x)^2+1=(x-x_1)^{m_1}\cdots(x-x_k)^{m_k}.
\end{align*}

Differentiating, we obtain $2P(x)P'(x)=(x-x_1)^{m_1-1} \cdots (x-x_k)^{m_k-1}Q(x)$, for some polynomial $Q(X)$. Hence $P'(x)$ is a multiple of $(x-x_1)^{m_1-1}\cdots(x-x_k)^{m_k-1}$. We conclude that 
$$m-1\geq \sum_i m_i-k=2m-k,$$ 
or, equivalently, $k\geq m+1$.

The degeneration in moduli occurs only when the geometric genus is $0$ or less. Since $\sum_i m_i=2m$, we can have 
$$\sum \lfloor \frac{m_i}{2} \rfloor \geq m-1$$
only if there are at most two odd numbers among $(m_1,\dots, m_n)$. As we are interested in points outside $\Delta_m$, there should be precisely two odd numbers, and the only case in which this occurs, under the additional restraint $k\geq m+1$, is when $$(m_1,\dots, m_n)=(\underbrace{2,\dots,2}_{m-1},1,1).$$ This happens only at the point $G=(\Delta_{m-1})_E\cap S_E$, which is unique by the proof of Lemma \ref{sectionG}.

\end{proof}

We would like now to understand the geometry of $\T_S\ra S_E$ around the points $(\Delta_{m-1})_E\cap S_E$ and $W_E\cap S_E$. We have the following 
\begin{lemma}
 Let $G=(\Delta_{m-1})_E\cap S_E$. Then $S_E$ has $m-1$ smooth branches around $G$.
\end{lemma}
\begin{proof} We work on $E_{(b_2)}\subset T_2$. The de-homogenization of the equation of $\T_{G}$ is 
$$y^2=\Psi_G(x)=P^2(x)+\lam_{2m}=(x-x_1)^2\cdots (x-x_{m-1})^2 S(x),$$
where $P(x)=x^m+\frac{1}{2}x^{m-2}+\sum_{i=3}^{m}\lam_ix^{m-i}$ is a polynomial of degree $m$ and $S(x)$ is a monic quadric. 

The tangent cone to $(\Delta_{m-1})_E$ at $G$ is a union of linear spaces of polynomials vanishing at the subset of $m-2$ nodes out of $x_1,\dots,x_{m-1}$. These linear spaces correspond to deformations preserving all nodes except one. We conclude that $(\Delta_{m-2})_E$ has $m-1$ smooth branches around $G$. To establish the lemma, it remains to show that $W_E$ intersects all branches transversely.

As we have seen in the proof of the Lemma \ref{sectionG}, the first order deformations of $P^2(x)+\lam_{2m}$ in $W_E$ are of the form 
$$\Psi_G(x)+\epsilon(2P(x)R(x)+\lam'_{2m})$$ where
$R(x)$ is a polynomial of degree $m-3$. 

It remains to observe that given a subset $\{x_{i_1},\dots,x_{i_{m-2}}\}$ of $\{x_1,\dots,x_{m-1}\}$, there is a unique, up to scaling, pair $(R(x), \lam'_{2m})$ such that
$$(x-x_{i_1})(x-x_{i_2})\cdots(x-x_{i_{m-2}}) \ \big{|} \ 2P(x)R(x)+\lam'_{2m}.$$
\end{proof}

It remains to discuss the geometry of $\T_S\ra S_E$ around a point $p \in \Delta\{2i, 2(m-i)\}_E$. The singularities of the fiber are tacnodes $y^2=x^{2i}$ and $y^2=x^{2(m-i)}$. We denote $T(0):=\Def(y^2=x^{2i})$ and $T(1):=\Def(y^2=x^{2(m-i)})$ to be the deformation spaces of these singularities. 
By Lemma \ref{versality}, the family $\T_S\ra S_E$ induces an isomorphism of a neighborhood of $p$ in $E$ with a neighborhood  of the origin in the product $T(0)\times T(1)$. Let $H_{0,2i}\subset T(0)$ and $H_{1,2(m-i)}\subset T(1)$ be the distinguished hyperplanes (see Section \ref{deformation space} for the definition). By Lemma \ref{versality-differential}, locally at $p$, variety $W_E$ is identified with a smooth $m-1$ dimensional subvariety $V$ of $T(0)\times T(1)$, and by Lemma \ref{versality-differential}, the tangent space of $V$ does not lie in the union of preimages of $H_{0,2i}$ and $H_{1,2(m-i)}$.

Therefore $S_E$ is identified with the curve $\Gamma_{i,m-i}$ residual to $\Delta_{i}\times \Delta_{m-i}$ 
in the intersection $$V\cap (\Delta_{i-1}\times \Delta_{m-i-1}).$$ This curve is an analog of a curve studied in the Lemma \ref{sectionG}.

Hence, to analyze the geometry of tails arising from a single tacnode, we are forced to consider the analogous problems posed for the product of the deformation spaces of multiple tacnodes. We do this in the next chapter. Note that to finish the analysis of $S$, we could simply use \cite[Lemma 4.3]{CH1} which describes the geometry of the curve $\Gamma_{i,m-i}$. However, we will need a slightly more general result for the case of several tacnodes, and so we reprove Lemma 4.3 of \cite{CH1} in our Lemma \ref{sectionGm}. Finally, we describe the geometry of the family $\T_S \ra S_E$ in Section \ref{resolving indeterminacy}.

\section{Multiple tacnodes case}\label{multiple tacnodes}

In this chapter, we study the geometry of the product of the deformation spaces of $n$ tacnodes.  Let $\m=(m_1,\dots,m_n)$ be a sequence of positive integer numbers. We define 
\begin{align*} 
M:=\prod_{i=1}^n m_i \ \ \text{and} \ \ m=\sum (m_i-1)+1.
\end{align*} 
For $1\leq i\leq n$, let $T(i)\cong \cc^{2m_i-1}$ be the versal deformation space of the tacnode $y^2=x^{2m_i}$, with coordinates
 $(a_{i, 2},\dots,a_{i, 2m_i})$. We set
$T=\prod_{i=1}^{n} T(i)$. Inside $T$, we denote the preimage of the hyperplane $a_{i,2m_i}=0$ by $H_i$.

In analogy with the case of a single tacnode, we first make a base change 
$$a_{i, j}=b_{i, j}^{j\frac{M}{m_i}}$$
to arrive at the space $T'=\spec \cc[b_{i,j}]$. 

We let $T''=Bl_{0} T'$, and denote by $E\subset T''$ the exceptional divisor of the blow-up.
In analogy with the single tacnode case, we also consider the quotients of $T''$ by groups 
$$\tilde{\mu}_{a,b}=\prod_{(i,j)\neq (a,b)} \mu_{j\frac{M}{m_i}};$$
we denote $T'' /\!/ \tilde{\mu}_{a,b}$ by $T_{a,b}$.


\begin{lemma}\label{tang_cone_2}
Consider any subvariety $W\subset T$, of dimension $m=\sum_{i=1}^n (m_i-1)+1$, that is smooth at the origin, contains $\Delta_\m$, and whose tangent space is not contained
 in the union of hyperplanes $H_i$. Then the exceptional divisor $W_E$ of $W$ in $T''$ is defined by equations
\begin{align}
a_{i,j}&=b_{i,j}^{j\frac{M}{m_i}}=0, \ \ \ 1\leq i \leq n, \ m_i+1\leq j \leq 2m_i-1; \label{eqt1}\\
a_{i,2m_i}&=\lam_i a_{1,2m_1}, \ \ \ 2\leq i\leq n; \label{eqt2}
\end{align}
where $\lam_i$ are some non-zero complex numbers. Here, $b_{i,j}$ are homogeneous coordinates on $E$, and $a_{i,j}=b_{i,j}^{j\frac{M}{m_i}}=0$.
\end{lemma}
\begin{proof}
The proof is a straightforward generalization of the proof of Lemma \ref{tang_cone}.
\end{proof}

\subsection{Geometry of $\Delta_{\m-\mathfrak{1}}$}
\begin{lemma}\label{sectionGm}
For any $W\subset T$, as in Lemma \ref{tang_cone_2}, we have
$$W\cap \Delta_{\m-\mathfrak{1}}=\Delta_{\m}\cup \Gamma_{\m},$$ 
where $\Gamma_{\m}$ is a curve intersecting $\Delta_\m$ with multiplicity $M$ at the origin.
\end{lemma}

\begin{proof}
First, we would like first to understand how $W_E$ and $(\Delta_{\m-\mathfrak{1}})_E$ intersect inside $E$. Any point of intersection outside $(\Delta_{\m})_E$ has homogeneous coordinates $[b_{i,j}]$ satisfying Equations \eqref{eqt1}-\eqref{eqt2}. Moreover, by Lemma \ref{poly}, $b_{i,2}\neq 0$. From now on, we work on the open affine 
\begin{align*}
E_{(b_{1,2})} &= \spec \cc[c_{i,j}]_{(i,j)\neq (1,2)}, \\
c_{i,j}&:=a^{i,j}/(b_{1,2})^{j\frac{M}{m_i}},
\end{align*}
inside the exceptional divisor of the map $f:T_{1,2}\ra T$. Note that $f$ is the weighted blow-up of $T$ followed by the base change of order $2m_2\cdots m_n$.  

Set theoretically, $W_E\cap (\Delta_{\m-\mathfrak{1}})_E$ consists of points whose coordinates satisfy the following conditions:

\begin{enumerate}
\item $c_{ij}$'s satisfy Equations \eqref{eqt1}-\eqref{eqt2}.
\item $c_{1,2m_1}\neq 0$. This simply means that a point is not in $\Delta_{\m}$.
\item The polynomial 
\begin{align*}
(x^{m_1}+x^{m_1-2}+c_{1,3}x^{m_1-3}\cdots+c_{1, m_1})^2+c_{1, 2m_1}
\end{align*} has $m_1-1$ double roots.
\item For each $2\leq i\leq n$, the polynomial    
\begin{align*}
(x^{m_1}+c_{i, 2}x^{m_1-2}+\cdots+c_{i, m_i})^2+c_{i, 2m_i}
\end{align*}
has $m_i-1$ double roots. 
\end{enumerate}

 The condition (3) together with Lemma \ref{poly} uniquely determine numbers $c_{1, j}$. Then Equation \eqref{eqt2} determines all numbers $c_{i, 2m_i}$, for $2\leq i\leq n$.
 Finally, we again invoke Lemma \ref{poly} and condition (4) to conclude that each of the vectors $(c_{i, 2}, c_{i, 3}, \dots, c_{i, 2m_i})$ is determined up to the scaling by an $m_i^{\th}$ root of unity. 
 
 Summarizing, $W_E\cap (\Delta_{\m-\mathfrak{1}})_E$ consists of $m_2\cdots m_n$ points. The tangent cone argument, analogous to that of the single tacnode case, shows that the intersection is transverse at these points. 

 We conclude that $\overline{\Gamma_{\m}}\cdot E=m_2\cdots m_n$. Observe that 
 \begin{align*}
 f^{*}\Delta_{\m}=\overline{\Delta_{\m}}+2m_1m_2\cdots m_n E .
 \end{align*}
 
 An application of the projection formula gives
 \begin{align*}
 \deg(f)\left(\Gamma_{\m}\cdot \Delta_{\m}\right) &=f(\overline{\Gamma_{\m}})\cdot \Delta_{\m} \\
 &=\overline{\Gamma_{\m}}\cdot f^{*}\Delta_{\m} \\ 
 &= \overline{\Gamma_{\m}}\cdot(\overline{\Delta_{\m}}+2m_1m_2\cdots m_nE) \\
 &=2m_1m_2\cdots m_n (\overline{\Gamma_{\m}}\cdot E) \\ 
 &=2m_1(m_2\cdots m_n)^2
 \end{align*}
 
 Since $\deg(f)=2m_2\cdots m_n$, we establish that
 \begin{align}\label{gm}
 \Gamma_{\m}\cdot \Delta_{\m}=\prod_{i=1}^{n} m_i=M.
 \end{align}
 
\end{proof}

We observe that a subvariety $W\subset T$ satisfying the conditions of Lemma \ref{tang_cone_2}, projects smoothly onto a subvariety $\text{proj}_k(W)$ of $T(k)$ satisfying the conditions of Lemma \ref{sectionG}. It makes sense, therefore, to 
talk about a curve $\Gamma_{k}$ inside $T(k)$, defined as the residual to $\Delta_{k,m_k}$ in the intersection $\text{proj}_k(W)\cap \Delta_{k,m_k-1}$.
We deduce the following
\begin{corollary}\label{projdegree} 
The projection $\text{proj}_k: \Gamma_{\m} \ra \Gamma_{k}$ has degree $\prod_{i\neq k}m_i$.
\end{corollary}

\begin{proof}
Clearly, $(\text{proj}_k)_{*} [\Gamma_{\m}]=\deg (\text{proj}_k)[\Gamma_k]$. Note that $(\text{proj}_k)^* \Delta_{i,m_i}=\Delta_{\m}$. Now using \eqref{gm} and applying the projection formula, we obtain the needed result.
\end{proof}

\subsection{The geometry of the miniversal family}\label{resolving indeterminacy}
In this section, we continue the study of the geometry of the surface $S^{d,\delta}(\al,\beta)$ and the family $\Y$ from Chapter \ref{local geometry}. We use the notations of the previous section and Section \ref{geometry indeterminacy} throughout.
We consider the intersection of a subvariety $W\subset T$, satisfying the conditions of Lemma \ref{tang_cone_2}, with $\Delta_{\m-\mathfrak{1}-e_1}$ inside $T$. We have 
$$W\cap \Delta_{\m-\mathfrak{1}-e_1}=\Delta_{\m} \cup S_1,$$
where $S_1$ is residual to $\Delta_{\m}$ in the intersection.  Generalizing the arguments of Chapter \ref{single tacnode}, we study the geometry of $S_1$ in what follows. 

Let $\Y_1 \subset \spec \cc[x,y]\times T''$ be the pullback to $T''$ of the miniversal family over $T(1)$. Recall that $\Y_1$ is 
given by the equation 
$$y^2=\Psi(a_{1,2},\dots, a_{1,2m_1})(x).$$

Define the ideal sheaf $\I=\biggl(\left(I_E^{M/m_1},x\right)^{m_1},y\biggr)$. Set 
$$\Z:=Bl_{\I} \Y \ \ \text{and} \ \ F: \Z\ra T''.$$

As in the previous chapter, $Bl_{\I} (\spec \cc[x,y]\times T'') \ra T''$ is a family of rational surfaces with fibers over $T'' \setminus E$ being affine planes, and fibers over the exceptional divisor $E$ being the union of $Bl_{(x^{m_1},y)}\cc^2$ and 
$\pp(1,\frac{M}{m_1},M)=\proj \cc[x,y,z]$. Here, $z$ stands for a local generator of $I_E$, and $\cc[x,y,z]$ is graded with $\deg x=\frac{M}{m_1}$, $\deg y=M$ and $\deg z=1$.

Consider a distinguished open affine of $T''$:
 \begin{align*}
 D_{(b_{1,2})}T' &'=\spec \cc[c_{i,j}], \ \text{where}  \\
 c_{1,2} &=b_{1,2}, \notag\\ 
 c_{i,j} &=\frac{b_{i,j}}{b_{1,2}}, \ (i,j)\neq (1,2). \notag
 \end{align*}
The exceptional divisor on $D_{(b_{1,2})}T''$ is given by $b_{1,2}=0$ and the restriction of $\Z$ to $E$ is 
$$\Z_E :=\Z \times_{T''} E=\mathcal{S}\cup \T ;$$
where $\mathcal{S}$ is
$$\{(y/x^{m_1})^2=1\}\subset E\times Bl_{(x^{m_1},y)}\cc^2$$ and $\T$ is given by
$$ \{y^2=\Psi\biggl(1,c_{1,3}^{3M/m_1},\dots, c_{1,2m_1}^{2M}\biggr) (x,z^{M/m_1}) \} \subset E\times \pp(1,M/m_1,M).$$

We denote the strict transform of $S_1$ in $T''$ by $\tilde{S}_1$ and the exceptional divisor of $\tilde{S}_1$ in $T''$ by $S_E$.
\begin{lemma} We have
$$S_E \cap (\Delta_{\m})_E=\bigcup_{a=1}^{\lfloor m_1/2 \rfloor} (\Delta\{2a,2(m_1-a)\}\times 0\times\cdots\times 0)_E.$$
\end{lemma}

\begin{proof}
This follows from the analogous statement for the single tacnode, proved in Lemma \ref{DeltamS}.
\end{proof}

\begin{remark} It follows that $S_E$ is a curve in $E$. Hence $S_1$ has pure dimension $2$.
\end{remark}

We now work on the quotient $T_{1,2}$ of $T''$ by $\tilde{\mu}_{(1,2)}$. Set $f: T_{1,2} \ra T''$ and $\Z_{1,2}=\Z /\!/\tilde{\mu}_{(1,2)}$. The distinguished open affine $D_{(b_{1,2})}T_{1,2}$ is isomorphic to $\spec[c_{i,j}]$, where
\begin{align*}
 c_{1,2}&=b_{1,2} \notag\\ 
 c_{i,j}&=\frac{a_{i,j}}{(b_{1,2})^{j\frac{M}{m_i}}}, \ (i,j)\neq (1,2). \notag
 \end{align*}
 
 The exceptional divisor $E_{(b_{1,2})}$ is given by $b_{1,2}=0$, and so is isomorphic to $\spec \cc[c_{i,j}]_{(i,j)\neq (1,2)}$. The quotient $\T_{1,2}$ of $\T$ is given by 
 equation 
 \begin{align*}
 y^2=\Psi(1,c_{1,3},\dots,c_{1,2m_1})(x,z).
 \end{align*}
 inside $\pp(1,\frac{M}{m_1},M)$.
 
 Note that at any point $p$ of $(\Delta\{2a,2(m_1-a)\}\times 0\times\cdots\times 0)_E$, the fiber of the family $\T_{1,2} \ra E_{(b_{1,2})}$ has two tacnodes: of order $a$ and $m_1-a$, respectively. 
 We denote their deformation spaces by $D(0)=\Def(y^2=x^{2a})$ and $D(1)=\Def(y^2=x^{2(m_1-a)})$. By a slight abuse of notation, we also denote $D(i)=\spec \cc[c_{i,j}]$, for $2\leq i \leq n$, and think of $D(i)$ as the deformation space of the order $m_i$ tacnode. 
 
 By Lemma \ref{versality}, the family $\T_{1,2}$ induces an isomorphism between a neighborhood of the point $p$ in $\spec \cc[c_{1,j}]_{3\leq j\leq 2m_1}$ and a neighborhood of the origin in $D(0)\times D(1)$. We conclude that at the point $p\in (\Delta\{2a,2(m_1-a)\}\times 0\times\cdots\times 0)_E$, 
there is a local isomorphism 
$$\Phi: E_{(b_{1,2})} \ra \prod_{i=0}^n D(i),$$
such that the map from $E_{(b_{1,2})}$ to the product of the versal deformation spaces of the singularities of $(\T_{1,2})_p$ is the composition of $\Phi$ and the projection to $D(0)\times D(1)$. By a slight reformulation of Lemma \ref{versality-differential}, the image $\Phi(W_E)$ satisfies the conditions of Lemma \ref{tang_cone_2}.

We can now state the result that we will use in the next chapter to compute $\discr{S_1}{\Delta_\m}{\lambda}$.
\begin{lemma}\label{curve identification}
Under the isomorphism $\Phi$, the exceptional divisor $S_E$ of $S_1$ is identified with the curve $\Gamma_{a,m_1-a,m_2,\dots,m_n}$, of Lemma \ref{sectionGm}, inside $\prod_{i=0}^n D(i)$.
\end{lemma}

Observe that the family $\Z$ over the blow-up $\tilde{S}_1$ the following property. Any $m^{\th}$-order tacnode, in any fiber, is a limit of $m-1$ nodes in the nearby fibers. Therefore, by the discussion in Chapter \ref{remarks}, the family $\Z \ra \tilde{S}_1$ has no points of indeterminacy along the preimage of $\Delta_{\m}$, at least after the normalization of the base. By abuse of notation, denote the strict transform of $S_1$ under $f$ also by $\tilde{S}_1$. It follows that the family $\Z_{1,2} \ra \tilde{S}_1$ has no point of indeterminacy along the preimage of $\Delta_{\m}$. 


\section{Calculation of the discrepancy}\label{calculation}

We continue the discussion of the previous chapter. Recall that we are considering the geometry of the branch $S$ of $S^{d,\delta}(\al,\beta)$ around the point of indeterminacy $[X]\in S^{d,\delta}(\al,\beta)\cap H_q$, and the geometry of the family $\Y \ra S$. The curve $X$ has $n$ ``new'' tacnodes $y_1,\dots, y_n$. The branch $S$ is defined as the closure of the deformations of $X$ such that $y_1$ deforms to $m_1-2$ nodes in the nearby fibers, and $y_i$ deforms to $m_i-1$ nodes, for $i\geq 2$. We use the notations of Section \ref{resolving indeterminacy} in what follows. 

In the neighborhood of the point $y_1$ on $X$, the family $\Y$ is the pullback of the family $\Y_1 \ra S_1$ under the isomorphism $\phi: S \ra S_1$. Recall, that we have constructed a blow-up
$f \co \tilde{S}_1 \ra S_1$, and a family $\Z_{1,2}\ra \tilde{S}_1$ with no points of indeterminacy along the preimage of $\Delta_{\m}$. We think of $f$ as the resolution of the moduli map along $\Delta_{\m}$. To compute the discrepancy $\discr{S_1}{\Delta_{\m}}{\lambda}$, we observe that 
$$f^*(\Delta_{\m}\cap S_1)=\overline{\Delta_{\m}}+2m_1m_2\cdots m_n S_E.$$
The coefficient of $S_E$ is $2m_1m_2\cdots m_n$ because $\Delta_{\m}$ is given by equation $a_{i,2m_i}=0$ (for any $i$) and the pullback of $a_{i,2m_i}=0$ to $T_{1,2}$ vanishes to the order $2m_1m_2\cdots m_n$ along $E$.

Recalling that $f$ is a map of degree $2m_2\cdots m_n$, we have 
\begin{align*} 
\discr{S_1}{\Delta_{\m}}{\lambda}&=\biggl((2m_1m_2\cdots m_n) S_E\cdot \lambda\biggr)/\deg f=m_1 S_E\cdot \lambda.
\end{align*}

It remains to compute $S_E\cdot \lambda$. Recall that the family $\Z_{1,2}$ over $S_E$ is a quotient of the union of two components. One component, $\mathcal{S}$, is an isotrivial family, and hence does not contribute to $S_E\cdot \lambda$. The other component is the family of tails $\T$, restricted to $S_E$. Its quotient is $\T_{1,2}$.

The family $\T_{1,2} \ra S_E$ is a family of generically hyperelliptic curves of arithmetic genus $m-1$, with generically $m-2$ nodes and two marked points. Consider the family 
$$\T_{\nu}:=\T_{1,2}\times_{S_E} (S_E)^{\nu} \ra (S_E)^{\nu}.$$ It induces a regular map from $(S_E)^{\nu}$ to $\Mg{1,2}$. By Proposition \ref{degenerate fibers}, the image curve intersects the boundary at the points of $(S_E)^{\nu}$ lying over
$$S_E\cap (\Delta_{\m-\mathfrak{1}})_E \ \ \text{and} \ \ S_E\cap (\Delta_{\m})_E.$$ 

At the points of $S_E\cap (\Delta_{\m-\mathfrak{1}})_E$, of which there are exactly $m_2\dots m_n$ by the proof of Lemma \ref{sectionGm}, the curve $S_E$ has $m_1-1$ smooth branches. Each branch intersects the boundary in $\Mg{1,2}$ transversely.

To calculate the multiplicity with which $S_E$ intersects the boundary at the points $\{p_1,\dots, p_{\kappa}\}$ of $(S_E)^{\nu}$ lying over the point 
$$p\in (\Delta\{2a,2(m_1-a)\}\times 0\times\cdots\times 0)_E,$$ we need to understand the geometry of the family $\T_{\nu}$ around the tacnodes in the fibers. By Lemma \ref{curve identification}, $S_E$ is identified, in a neighborhood of $p$, with the curve $\Gamma_{a,m_1-a,m_2,\dots,m_n}$ inside $\prod_{i=0}^n D(i)$. Therefore, by Corollary \ref{projdegree}, the curve $(S_E)^{\nu}$ maps with degree $(m_1-a)m_2\dots m_n$ onto the curve $\Gamma_a$ inside the deformation space 
$D(0)$ of the tacnode $y^2=x^{2a}$. Let $r_i$, for $1\leq i\leq \kappa$, be the ramification index of the map
$$(S_E)^{\nu}\ra \Gamma_a$$
at the point $p_i$.
Then by Lemma \ref{delta-multiplicity}, the surface $\T_{\nu}$ has an $A_{r_i-1}$-singularity at the $a^{\th}$-order tacnode in the fiber $(\T_{\nu})_{p_i}$. Hence, this $a^{\th}$-order tacnode contributes $r_i$ to the intersection multiplicity of $(S_E)^{\nu}$ with $\Delta$ at the point $p_i$. Remembering that
 $$\sum_{i=1}^{\kappa} r_i=(m_1-a)m_2\dots m_n,$$ 
and applying the same argument to the other tacnode, we conclude that the intersection number of $(S_E)^{\nu}$ with the boundary at points lying over $p$ is
 $$(m_1-a)m_2\cdots m_n+am_2\cdots m_n=\prod_{i=1}^n m_i.$$
 
Noting that $S_E\cap (\Delta_{\m})_E$ has exactly $m_1-1$ points, we sum up the contributions of $S_E\cap (\Delta_{\m-\mathfrak{1}})_E$ and $S_E\cap (\Delta_{\m})_E$ to the intersection number of $(S_E)^{\nu}$ with the boundary in $\Mg{1,2}$:
\begin{align*}
(S_E)^{\nu}\cdot \Delta&=(m_1-1)(m_1m_2\cdots m_n)+(m_2\cdots m_n)(m_1-1) \\
&=(m_1^2-1)(m_2\cdots m_n). 
\end{align*} 

Therefore, 
\begin{align*} 
\discr{S_1}{\Delta_{\m}}{\lambda}&=m_1S_E\cdot \lambda = m_1 ((S_E)^{\nu}\cdot \Delta)/12 \\
&=(m_1^2-1)(m_1m_2\cdots m_n)/12 .
\end{align*}

Performing the same calculation for the surfaces $S_i$, residual to $\Delta_{\m}$ in the intersection 
$\Delta_{\m-\mathfrak{1}-e_i}\cap W,$
we arrive at the formula
 $$\discr{S_\m}{\Delta_\m}{\lambda}=\sum_{i=1}^n \discr{S_i}{\Delta_{\m}}{\lambda}=\frac{1}{12}\left(\prod_{i=1}^n m_i\right)\left(\sum_{i=1}^n (m_i^2-1)\right).$$

Together with Equations \eqref{discrepancy-total} and \eqref{L-induction2}, this finishes the proof of Theorem \ref{main theorem}.

\section{Examples and applications}

\subsection{Slopes of effective divisor on $\Mg{g}$.} Define a slope of a curve $C\subset \Mg{g}$ by 
$$s(C):=\frac{C \cdot \Delta}{C\cdot \lambda} \ .$$ 
If $C$ is a moving curve in $\Mg{g}$, then the slope of $C$ gives a lower bound on the slope $s_g$ of effective divisors on $\Mg{g}$ as defined in \cite{HM}.   
We have noted in the introduction that $C^{d,\delta}_{\irr}$ is a moving curve in $\Mg{g}$ when $3d-2g-6\geq 0$.

We can slightly relax condition $3d-2g-6\geq 0$ by considering curves $C_{\irr}^{d,\delta}$ with $3d-2g-6=-1$. By irreducibility of $V^{d,\delta}$ (see \cite{Harris1}) and by work of Eisenbud and Harris (see \cite{EH}), we know that, in this case, deformations of $C^{d,\delta}_{\irr}$ span an irreducible Brill-Noether divisor inside $\Mg{g}$ whose slope is $6+\frac{12}{g+1}$. Therefore, by a standard argument (cf. \cite{HM}), whenever $3d-2g-6\geq -1$, we have a bound:
$$s_g\geq \min \{6+\frac{12}{g+1}, s(C^{d,\delta}_{\irr})\},$$
here, as always, $\delta=\binom{d-1}{2}-g$.

Using the recursion of Theorem \ref{main theorem}, we computed the slopes of curves $C_{\irr}^{d,\delta}$ for $d\leq 16$ with a help of a computer. This gives us some lower bounds on $s_g$ for $g\leq 21$. 
\begin{table}[htdp]
\caption{Slopes of $C^{d,\delta}_{\irr}$}
\begin{center}
\begin{tabular}{|c|c|c|c|c||c|c|c|c|c|}
\hline
$g$ & $d$ & $\delta$  & $3d-2g-6$ & $s(C^{d,\delta}_{\irr})\approx$ & $g$ & $d$  & $\delta$  & $3d-2g-6$ & $s(C^{d,\delta}_{\irr})\approx$ \\
\hline
2 & 4 & 0 & 1 &  10               &               12 & 10 & 24 & 0 & 6.29              \\ 	
\hline
3 & 4 & 0 & 0 & 9              &    		         13 & 11 & 32 & 1 & 6.03 \\
\hline
4 & 5 & 2 & 1 & 8.45              & 			14 & 11 & 31 & -1 & 6.00 \\	
\hline
5 & 5 & 1 & -1 & 8.16              & 			       15 & 12 & 40 & 0 & 5.76 \\	
\hline
6 & 6 & 4 & 0 & 7.76              & 				16 & 13 & 50 & 1 & 5.53 \\
\hline
7 & 7 & 8 & 1 & 7.37              & 				17 & 13 & 49 & -1 & 5.52 \\
\hline
8 & 7 & 7 & -1 & 7.27              & 				18 & 14 & 60 & 0 & 5.31 \\
\hline
9 & 8 & 12 & 0 & 6.93              & 				19 & 15 & 72 & 1 & 5.12 \\
\hline
10 & 9 & 18 & 1 & 6.62              & 				20 & 15 & 71 & -1 & 5.11 \\
\hline
11 & 9 & 17 & -1 & 6.57              & 				21 & 16 & 84 & 0 & 4.93 \\
\hline

\end{tabular}
\end{center}
\label{default}
\end{table}

Note that the bounds for $g=2$ and $g=3$ are sharp. The bounds for $g=4$ and $g=5$ are better than those given in \cite{HM}, but are still not sharp. 
We remark that there are examples of moving curves in $\Mg{g}$, providing sharp lower bounds, at least for small $g$. For $g\leq 6$, see, for example, \cite{CHS}. Finally, even though we have nothing to say about the asymptotic behavior of the bounds produced by curves $C^{d,\delta}_{\irr}$, it would not be surprising if these bounds approached $0$, as $g$ approached $\infty$.

\subsection{Codimension one numbers on $V_{d,g}$}

 Consider the Severi variety $V_{d,g}$ of irreducible plane curves of degree $d$ and geometric genus $g$. In \cite{DH1}, Diaz and Harris have computed a great number of geometrically meaningful divisors on $V_{d,g}$ in terms of three standard classes (see loc. cit. for the notations):
 \begin{align*} A&=\pi_*(\omega^2), && B=\pi_*(\omega\cdot D), && C=\pi_*(D^2),
 \end{align*}
 and boundary divisors $\Delta_0$ and $\Delta_{i,j}$. 
 In particular, we recall the formulas for the classes of divisors of curves with cusps, triple points and tacnodes:
\begin{align}
CU &= 3A+3B+C-\Delta, \label{cusps}\\
TN &= (3(d-3)+2g-2)A+(d-9)B-\frac{5}{2}C+\frac{3}{2}\Delta, \label{tacnodes}\\
TR &= \biggl(\frac{d^2-6d+8}{2}-g+1\biggr)A-\frac{d-6}{2}B+\frac{2}{3}C-\frac{1}{3}\Delta. \label{triple}
\end{align}

The number of curves in $C^{d,\delta}_{\irr}$ with any codimension one behavior, which was studied in \cite{DH1}, is
 expressed in terms of intersection numbers of $C^{d,\delta}_{\irr}$ with $A,B,C$ and $\Delta$. 
 
  We recall that $C^{d,\delta}\cdot A=N^{d,\delta}$, and so can be computed by Caporaso's and Harris's recursion. The number $C^{d,\delta}\cdot B$ can also be expressed in terms of 
 degrees of Severi varieties (see \cite{RV-en}). However, intersection numbers $C^{d,\delta}\cdot C$ were not known. 
 
 By Mumford's formula $$C=12\lambda-\Delta$$ and hence 
 the intersection $C^{d,\delta}\cdot C$ is computed in terms of $L^{d,\delta}=C^{d,\delta}\cdot \lambda$ and intersections of $C^{d,\delta}$ with boundary divisors.
 
We remark that numbers  
\begin{align} B^{d,\delta}(\al,\beta):=C^{d,\delta}(\al,\beta)\cdot B\end{align}
can be computed using the following recursion.
\begin{proposition}\label{B-recursion} The numbers $B^{d,\delta}(\alpha,\beta)$ satisfy the following recursion:
\begin{align*} B^{d,\delta}(\alpha,\beta) &= \sum_k kB^{d,\delta}(\alpha+e_k,\beta-e_k) \\
&+\sum  \bigg( I^{\beta'-\beta}\binom{\alpha}{\alpha'}\binom{\beta'}{\beta}B^{d-1,\delta'}(\alpha',\beta') \\
&+2\sum_k I^{\beta'-\beta}\binom{\alpha}{\alpha'}\binom{\beta'-e_k}{\beta}
N^{d-1,\delta'}(\alpha'+e_k,\beta'-e_k) \bigg),
\end{align*}
where the second sum is taken over all triples $(\delta',\alpha',\beta')$ satisfying $|\beta'-\beta|+\delta-\delta'=d-1$.
\end{proposition}
\begin{proof}
The formula follows from definitions and Theorem \ref{CH hyperplane section theorem}.
\end{proof}

Theorems \ref{CH hyperplane section theorem} and \ref{main theorem}, together with Proposition \ref{B-recursion}, allow us to compute inductively the intersection numbers of 
$C^{d,\delta}_{\irr}$ with all the standard divisor classes and, hence, to find solutions to a large class of codimension one enumerative problems on $V^{d,\delta}_{\irr}$. For
example, one can compute the number of irreducible plane curves of degree $d$ with either a single cusp, a tacnode, or a triple point, and nodes as the only other singularities, passing through the appropriate number of general points in $\pp^2$. 

\subsection{Enumerative applications}
\begin{example}
A simple induction argument involving the recursion of Theorem \ref{main theorem} shows that 
\begin{align*}
L^{d,1} &= \frac{3}{2}(d-1)(d-2)(d-3)(d+1), \\
L^{d,2} &=\frac{1}{4}(9d^6-63d^5+66d^4+333d^3-553d^2-480d+828).
\end{align*}
We skip the proof of these formulas, and only note that it requires an introduction of auxiliary functions, such as $N^{d,0}(0,(d-4,2))$ and $L^{d,0}(0,(d-2,1))$. The derivation of the closed form formula for these functions is inductive as well. 

Similarly, we have 
\begin{align*}
B^{d,1} &=3(d-1)(2d^2-5d+1)\label{degB}, \\
N^{d,2} &= \frac{1}{2}(9d^4-36d^3+12d^2+81d-33).
\end{align*}
Using Equation \eqref{cusps} and the equality $\Delta\cdot C^{d,1}=2N^{d,2}$, we recover the number of cuspidal curves in $C^{d,1}$:
\begin{align*}
CU=3N^{d,1}+3B^{d,1}+12L^{d,1}-4N^{d,2}=12(d-1)(d-2).
\end{align*}
The right-hand side is the classical formula for the degree of the cuspidal locus in $\pp(d)$.
\end{example}

The numerical computations suggest that, for a fixed $\delta$, the function $L^{d,\delta}$ is a polynomial of degree $2\delta+2$ in $d$, with the leading term $\frac{3^\delta}{2\delta!}$. 
This should be compared with the G\"{o}ttsche's conjecture (\cite{Gott}) that states that $N^{d,\delta}$ is a polynomial of degree $2\delta$ in $d$. The conjecture was proved for $\delta\leq 8$ by Kleiman and Piene in \cite{KP-np}. 

In \cite[Theorem (1.2)]{KP-sing}, Kleiman and Piene enumerate curves  with certain singularities in a sufficiently ample linear series on an arbitrary surface. In particular, their formulas compute the number of degree $d$ plane curves with a single triple point and $\delta$ nodes, where $\delta\leq 3$, passing through an appropriate number of general points. The postulated number is a polynomial of degree $2(\delta+1)$ in $d$. Our methods allow us to compute the number of plane curves with a triple point and an arbitrary number of nodes. However, we cannot produce a closed form formula. We ran computations, with a help of a computer, for $d\leq 13$ and $\delta\leq 3$, and our numbers agree with those of \cite{KP-sing}. 

In a different direction, for $g=0,1,2,3$, the recursions for the codimension one ``characteristic'' numbers of plane curves were given by Vakil in \cite{RV-en}. We note that our computations agree with those of Vakil, presented in the table on page 19 of loc. cit. 

\bibliography{Bib}
\bibliographystyle{alpha}

\end{document}